\documentclass[a4paper,11pt]{amsart}
\usepackage{amsmath,amsthm,amsfonts,amssymb}

\usepackage{ifpdf} 
\ifpdf       
\usepackage[pdftex, linktocpage]{hyperref}
\else
\usepackage[hypertex, linktocpage]{hyperref}
\fi
\usepackage{verbatim}

\title{Infinite paths and cliques in random graphs} 
\author{Alessandro Berarducci, Pietro Majer \and Matteo Novaga}
\address{Alessandro Berarducci \\
        Dipartimento di Matematica, Universit\`a di Pisa,
        Largo B. Pontecorvo 5, 56127 Pisa, Italy,
        email: \textit{berardu@dm.unipi.it}}
\address{Pietro Majer \\
        Dipartimento di Matematica, Universit\`a di Pisa,
        Largo B. Pontecorvo 5, 56127 Pisa, Italy, 
        email: \textit{majer@dm.unipi.it}}
\address{Matteo Novaga \\
        Dipartimento di Matematica, Universit\`a di Padova,
        Via Trieste 63, 35121 Padova, Italy, 
        email: \textit{novaga@math.unipd.it}}
 
\date{25 Mar. 2011}       

\DeclareMathOperator{\R}{\mathbb R}
\DeclareMathOperator{\NN}{\mathbb N}

\DeclareMathOperator{\im}{Im}

\DeclareMathOperator{\X}{\mathbb X} 
\newcommand{\comb}[2]{#1^{(#2)}}

\newcommand{\st}{\, : \,}

\newcommand{\A}{\mathcal A}
\newcommand{\Sc}{\mathfrak S_c(\NN)}
\newcommand{\Ij}{{\rm Inj}(\NN)}
\newcommand{\Ic}{{\rm Incr}(\NN)}

\newcommand{\nin}{\not\in}
\newcommand{\eps}{\varepsilon}
\newcommand{\ca}[1]{{\mathcal #1}}
\newcommand{\rest}[1]{\hspace {-0.2em} \upharpoonright_{#1}}

\theoremstyle{plain}
\newtheorem{theorem}{Theorem}
\newtheorem{lemma}[theorem]{Lemma}
\newtheorem{proposition}[theorem]{Proposition}
\newtheorem{corollary}[theorem]{Corollary}
\newtheorem{conjecture}[theorem]{Conjecture}

\newtheorem{question}[theorem]{Question}
\newtheorem{fact}[theorem]{Fact}
\theoremstyle{definition}
\newtheorem{remark}[theorem]{Remark}
\newtheorem{definition}[theorem]{Definition}
\newtheorem{example}[theorem]{Example}
\newtheorem{exercise}[theorem]{Exercise}
\newtheorem{problem}{Problem}

\numberwithin{theorem}{section}
\numberwithin{equation}{section}

\newcommand{\bt}{\begin{theorem}}
\newcommand{\et}{\end{theorem}}

\newcommand{\bl}{\begin{lemma}}
\newcommand{\el}{\end{lemma}}

\newcommand{\bd}{\begin{definition}}
\newcommand{\ed}{\end{definition}}

\newcommand{\beq}{\begin{equation}}
\newcommand{\eeq}{\end{equation}}

\newcommand{\bexa}{\begin{example}}
\newcommand{\eexa}{\end{example}}

\newcommand{\bexe}{\begin{exercise}}
\newcommand{\eexe}{\end{exercise}}

\newcommand{\bfact}{\begin{fact}}
\newcommand{\efact}{\end{fact}}

\newcommand{\bprop}{\begin{proposition}}
\newcommand{\eprop}{\end{proposition}}

\newcommand{\bp}{\begin{proof}}
\newcommand{\ep}{\end{proof}}

\newcommand{\bc}{\begin{corollary}}
\newcommand{\ec}{\end{corollary}}

\newcommand{\bq}{\begin{question}}
\newcommand{\eq}{\end{question}}

\newcommand{\bconj}{\begin{conjecture}}
\newcommand{\econj}{\end{conjecture}}

\newcommand{\br}{\begin{remark}}
\newcommand{\er}{\end{remark}}

\newcommand{\bproblem}{\begin{problem}}
\newcommand{\eproblem}{\end{problem}}

\newcommand{\ben}{\begin{enumerate}}
\newcommand{\een}{\end{enumerate}}

\newcommand{\p}{\varphi}
\renewcommand{\phi}{\varphi}

\begin{document}

\begin{abstract} We study the thresholds for the emergence of various properties in random subgraphs of $(\NN, <)$. In particular, we give sharp sufficient conditions for the existence of (finite or infinite) cliques and paths in a random subgraph. 
No specific assumption on the probability is made. 
The main tools are a topological version of Ramsey theory, exchangeability theory and elementary ergodic theory. 
\end{abstract}

\maketitle

\tableofcontents

\section{Introduction}
In this paper we introduce a new method in order to deal with some combinatorial problems in random graphs, originally proposed in \cite{EH:64}. 
Some of this questions have been successfully addressed in \cite{FT:85}, using different techniques.  
We obtain new and self-contained proofs of some of the results in \cite{FT:85}; moreover 
with this method we expect to be able to treat similar problems in more general random graphs.

Let $G=(\NN, \comb {\NN} 2)$ be the directed graph over $\NN$ with set of edges $\comb {\NN} 2:=\{(i,j) \in \NN^2 \st i<j\}$. Let us randomly choose some of the edges of $G$, that is, we associate to the edge $(i,j)\in\comb {\NN} 2$ a measurable set $\X_{i,j}\subseteq \Omega$, where $(\Omega, \mathcal A, \mu)$ is a base probability space. Assuming $\mu(\X_{i,j})\geq \lambda$ for each $(i,j)$, we then ask whether the resulting random subgraph $\X$ of $(\NN, \comb {\NN} 2)$ contains an infinite path:

\bproblem \label{increasing} Let $(\Omega, \mathcal A, \mu)$ be a probability space. Let $\lambda > 0$ and for all $(i,j) \in \comb {\NN} 2$, let $\X_{i,j}$ be a measurable subset of $\Omega$ with $\mu(\X_{i,j}) \geq \lambda$. 
Is there an infinite increasing sequence $\{n_i\}_{i\in\NN}$ 
such that $\bigcap_{i\in \NN} \X_{n_i, n_{i+1}}$ is non-empty?  
\eproblem

More formally, a random subgraph $\X$ of a directed graph $G= (V_G, E_G)$ (with set of edges $E_G \subset V_G\times V_G$), is a measurable function $\X: \Omega \to 2^{E_G}$
where $\Omega = (\Omega, \ca A, \mu)$ is a probability space, 
and $2^{E_G}$ is the powerset of $E_G$, identified with the set of all functions from $E_G$ to $\{0,1\}$ (with the product topology and the $\sigma$-algebra of its Borel sets). For each $x\in \Omega$, we identify $\X(x)$ with the subgraph of $G$ with vertices $V_G$ and edges $\X(x)$. Given $e\in E_G$, the set $\X_e:= \{x\in \Omega \st e \in \X(x)\}$ represents the event that the random graph $\X$ contains the edge $e\in E_G$.  
The family $(\X_e)_{e\in E_G}$ determines $\X$ putting: $\X(x) = \{e \in E_G \st x \in \X_e\}$. So a random subgraph of $G$ can be equivalently defined as a function from $E_G$ to $2^{\Omega}$ assigning to each $e\in E_G$ a measurable subset $\X_e$ of $\Omega$.

As in classic percolation theory, we wish to estimate the probability that $\X$ contains an infinite path, in terms of a parameter $\lambda$ that 
bounds from below the probability $\mu(\X_e)$ that an edge $e$ belongs to $\X$. 
Note that it is not a priori obvious that the existence of an infinite path has a well-defined probability, since it corresponds to the uncountable union of the sets $\bigcap_{k\in\NN}\X_{i_k, i_{k+1}}$ 
over all strictly increasing sequences $i:\,\NN\to\NN$. 
However, it turns out that it belongs to the $\mu$-completion
of the $\sigma$-algebra generated by the $\X_{i,j}$. 
It has to be noticed that the analogy with classic bond percolation 
is only formal, the main difference being that in the usual percolation models (see for instance \cite{G:99}) the events $\X_{i,j}$ are supposed 
\emph{independent}, whereas in the present case the probability distribution is completely general, i.e. we do not impose any restriction on the events $\X_{i,j}$, and on the probability space $\Omega$.   

Problem \ref{increasing} has been originally
proposed by P. Erd\H os and A. Hajnal in \cite{EH:64}, and an answer was given  
by D. H. Fremlin and M. Talagrand in \cite{FT:85}, 
where other related and more general problems are also considered. 
In particular they show that the threshold for the existence of infinite paths is $\lambda=1/2$, 
under the assumption that the probability space $(\Omega, \mathcal A, \mu)$ is $[0,1]$ equipped with the Lebesgue measure (although the extension to a general probability space should not be difficult). 
One of the main goals of this paper is to present a  general method, different from the one in \cite{FT:85}, 
which in particular allows us to recover the same result as in \cite{FT:85} (see Theorem \ref{infpath}).
Our approach relies on the reduction to the following dual problem:

\bproblem \label{dual}
Given a directed graph $F$, determine the minimal
$\lambda_c$ such that, whenever $\inf_{e\in \comb {\NN} 2} \mu(\X_e) > \lambda_c$, there is a graph morphism $f\colon \X(x)\to F$ for some $x\in \Omega$. 
\eproblem 

Problem \ref{increasing} can be reformulated in this setting by letting
$F$ be the graph $(\omega_1, >)$ where $\omega_1$ is the first uncountable ordinal. This depends on the fact that a subgraph $H$ of $(\NN,\comb {\NN} 2)$ does not contain an infinite path if and only if it admits a rank function with values in $\omega_1$. 
Therefore, if a random subgraph $\X$ of $(\NN, \comb {\NN} 2)$ has no infinite paths, it is defined a $\mu$-measurable map $\varphi \colon \Omega \to \omega_1^{\NN}$ where $\varphi(x)(i)$ is the rank of the vertex $i\in \NN$ in the graph $\X(x)$. 
It turns out that $\phi_\#(\mu)$ is a compactly supported Borel measure on $\omega_1^{\NN}$,
and that $\phi(\X_{i,j})\subseteq A_{i,j}:= 
\left\{ x\in \omega_1^{\NN}:\, x_i>x_j\right\}$. 
As a consequence, in the determination of the threshold for existence of infinite paths 
\begin{equation}\label{eqlac}
\lambda_c:=\sup\left\{\inf_{(i,j)\in \comb {\NN} 2}\mu(\X_{i,j}):\ \X \ \textrm{random graph without infinite paths}\right\}
\end{equation}
we can set $\Omega=\omega_1^{\NN}$, $\X_{i,j}=A_{i,j}$, and reduce to 
the variational problem on the convex set ${\mathcal M}^1_c(\omega_1^{\NN})$ of compactly supported probability measures on $\omega_1^{\NN}$:
\begin{equation}\label{probvar}
\lambda_c=\sup_{m\in {\mathcal M}^1_c(\omega_1^{\NN})}\inf_{(i,j)\in \comb {\NN} 2}\,m\left( A_{i,j}\right).
\end{equation}
As a next step, we show that in \eqref{probvar} we can equivalently take the supremum in the smaller class of all the compactly supported \emph{exchangeable measures} on $\omega_1^{\NN}$ 
(see Appendix \ref{app-exch} and references therein for a precise definition).
Thanks to this reduction, we can explicitly compute $\lambda_c=1/2$ ( Theorem \ref{infpath}).
We note that the supremum in \eqref{probvar} is not attained, which implies that for $\mu(\X_{i,j})\ge 1/2$ infinite paths occurs with positive probability. 

In Section \ref{secrelp}, we consider  again Problem \ref{dual} and we give a complete solution when $F$ is a finite graph, showing in particular that 
\[
\lambda_c = \sup_{\lambda\in \Sigma_F}\sum_{(a,b)\in E_F}\lambda_a\lambda_b
\]
where $\Sigma_{F}$ is the set of all sequences $\{\lambda_a\}_{a\in V_F}$ with values in 
$[0,1]$ and such that $\sum_{a\in V_F}\lambda_a=1$. By the appropriate choice of $F$ we can determine the threholds for the existence of paths of a given finite length   (Section 3 and Remark \ref{implicit}), or for the property of having chromatic number $\geq n$ (Section 6). 

We can consider Problems \ref{increasing} and \ref{dual} for a random subgraph $\X$ of an arbitrary directed graph $G$, not necessarily equal to $(\NN, \comb {\NN} 2)$. 
However, it can be shown that, if we replace $(\NN, \comb {\NN} 2)$ with a 
finitely branching graph $G$ (such as a finite dimensional network), the probability that $\X$ has an infinite path may be zero even if $\inf_{e\in E_G} \mu(\X_e)$ is arbitrarily close to $1$ ( Proposition \ref{fin-b}). Another variant is to consider subgraphs of $\comb {\R} 2$ rathen than $\comb {\NN} 2$ but it turns out that this makes no difference in terms of the threshold for having infinite paths in random subgraphs (Remark \ref{reals}). 

In Section \ref{secproblem} we fix again $G=(\NN,\comb {\NN} 2)$ and we ask if a random subgraph $\X$ of $G$ contains an infinite clique, i.e. a copy of $G$ itself. More generally we consider the following problem. 

\bproblem \label{problemstrong}
Let $(\Omega, \A,\mu)$ be a probability space. Let $\lambda > 0$ and for all $(i_1,\ldots,i_k) \in \comb {\NN} k$, 
let $\X_{i_1, \ldots, i_k}$ be a measurable subset of $X$ with $\mu(\X_{i_1, \ldots, i_k})\geq \lambda$. 
Is there an infinite set $J\subset\NN$ such that $\bigcap_{(i_1,\ldots,i_k)\in \comb J k} \X_{i_1, \ldots, i_k}$ is non-empty? 
\eproblem

This problem is a random version of the classical Ramsey theorem \cite{Ra} (we refer to \cite{GP:73,PR:05}, and references therein, for various generalization of Ramsey theorem). Clearly Ramsey theorem implies that the answer to Problem \ref{problemstrong} is positive when $\Omega$ is finite. Moreover it can be shown that the answer remains positive when $\Omega$ is countable (Example \ref{excountable}). However when $\Omega = [0,1]$ (with the Lebesgue measure) the probability that $\X$ contains an infinite clique may be zero even when $\inf_{e\in E_G} \mu(\X_e)$ is arbitrarily close to $1$ (see Example \ref{exo}). We will show that Problem \ref{problemstrong} has a positive answer if the indicator functions of the sets 
$\X_{i_1, \ldots, i_k}$ all belong to a compact subset of $L^1(\Omega, \mu)$ (see Theorem \ref{taiut}).

 Our original motivation for the above problems came from the following situation.
Suppose we are given a space $E$ and a certain family $\Omega$ of sequences on $E$ 
(e.g., minimizing sequences of a functional, or orbits of a discrete dynamical system, etc). A typical general problem asks for existence of 
a sequence in the family $\Omega$, that admits a subsequence with a prescribed property. One approach to it is by means of measure theory.
The archetypal situation here come from recurrence theorems: one may ask if there exists a subsequence
which belongs frequently to a given subset $C$ of the ``phase'' space $\Omega$ (we refer to such sequences as ``$C$-recurrent orbits'').
If we consider the set $\X_i: =  \{x \in  \Omega : x_i \in C\}$, 
then a standard sufficient condition for existence of
$C$-recurrent orbits is $\mu(\X_i) \geq \lambda >0$, 
for some probability measure $\mu$ on $\Omega$. 
In fact is easy to check that the set of $C$-recurrent orbits has measure at least $\lambda$ 
by an elementary version of a Borel-Cantelli lemma
(see Proposition \ref{lemone}). This is indeed the existence argument in the Poincar\'e Recurrence Theorem for measure preserving transformations. 
A more subtle question arises when one looks for 
a subsequence satisfying a given relation between two successive (or possibly more) terms: given a subset $R$ of $E \times E$
we look for a subsequence $x_{i_k}$ such that 
$(x_{i_k}, x_{i_{k+1}}) \in R$ for all $k\in\NN$. As before, we may
consider the subset of $\Omega$, with double indices $i<j$,
$\X_{i,j} := \{ x\in \Omega: (x_i,x_j) \in R\}$ and we are then led to Problem \ref{increasing}.

\section{Notations} \label{notations}

We follow the set-theoretical convention of identifying a natural number $p$ with the set $\{0,1,\ldots, p-1\}$ of its predecessors. More generally an ordinal number $\alpha$ coincides with the set of its predecessors. With these conventions the set of natural numbers $\NN$ coincides with the least infinite ordinal $\omega$. As usual $\omega_1$ denotes the first uncountable ordinal, namely the set of all countable ordinals. 

Given two sets $X,Y$ we denote by $X^Y$ the set of all functions from $Y$ to $X$. If $X,Y$ are linearly ordered we denote by $\comb X Y$ the set of all increasing functions from $Y$ to $X$. In particular $\comb {\NN} p$ (with $p\in \NN$) is the set of all increasing $p$-tuples from $\NN$, where a $p$-tuple $\boldsymbol i = (i_0,\ldots, i_{p-1})$ is a function $\boldsymbol i \colon p \to \NN$. The case $p=2$, with the obvious identifications, takes the form $\comb {\NN} 2 = \{(i,j)\in \NN^2 \st i<j\}$. 

Any function $f\colon X \to X$ induces a function $f_* \colon X^Y \to X^Y$ by $f(u) = f \circ u$. On the other hand a function $f \colon Y \to Z$ induces a function $f^* \colon X^{Z} \to X^Y$ by $f^* (u) = u \circ f$. In particular if ${\sf S} \colon \NN \to \NN$ is the successor function, ${\sf S}^*  \colon X^{\NN} \to X^{\NN}$ is the {\em shift map}. 

We let $\Sc, \Ij, \Ic\subset \NN^{\NN}$ be the families of maps $\sigma: \NN\to \NN$
which are compactly supported permutations\footnote{that is, finite perturbations of the identity}, injective functions and strictly increasing functions, respectively. Note that with the above conventions $\Ic = \comb {\NN} \omega$. 

Given a measurable function $\psi \colon X \to Y$ between two measurable spaces and given a measure $m$ on $X$, we denote as usual by $\psi_\#(m)$ the induced measure on $Y$. 

Given a compact metric space $\Lambda$, the space $\mathcal M(\Lambda^{\NN})$ of Borel measures on $\Lambda^{\NN}$ can be identified with $C(\Lambda^{\NN})^*$, i.e. the dual of the Banach space of all continuous functions on $\Lambda^{\NN}$. 
By the Banach-Alaoglu theorem the subset ${\mathcal M}^1(\Lambda^{\NN})\subset \mathcal M(\Lambda^{\NN})$ of probability measures is a compact (metrizable) subspace of $C(\Lambda^{\NN})^*$ endowed with the weak$^*$ topology.

Given $\sigma \colon \NN \to \NN$ we have $\sigma^* \colon \Lambda^{\NN} \to \Lambda^{\NN}$ and $\sigma^*_\# \colon 
{\ca M}^1(\Lambda^{\NN}) \to {\ca M}^1(\Lambda^{\NN})$.  To simplify notations we also write $\sigma \cdot m$ for $\sigma^*_\# m$. Note the contravariance of this action: 
\begin{equation} 
\theta \cdot \sigma \cdot m = (\sigma \circ \theta)\cdot m  \,.
\end{equation}
Similarly given $r\in \NN$ and $\iota \in \comb {\NN} r$, we have 
$\iota^*_\# \colon 
{\ca M}^1(\Lambda^{\NN}) \to {\ca M}^1(\Lambda^{r})$ and 
we define $\iota \cdot m = \iota^*_\#(m)$. 

Given a family $\ca F \subset \NN^{\NN}$, we say that $m$ is $\ca F$-invariant if $\sigma \cdot m = m$ for all $\sigma \in \ca F$. 

\section{Finite paths in random subgraphs} 


As a preparation for the study of infinite paths (Problem \ref{increasing}) we first consider the case of finite paths.
The following example shows that there are random subgraphs $\X$ of $(\NN, \comb {\NN} 2)$ such that $\inf_{e\in \comb {\NN} 2} \X_e$ is arbitrarily close to $1/2$, and yet $\X$ has probability zero of having infinite paths. 

\bexa\label{ex4}  Let $p\in \NN$ and let $\Omega = p^{\NN}$ with the Bernoulli probability measure $\mu = B_{(1/p,\ldots,1/p)}$. For $i<j$ in $\NN$ let 
$\X_{i,j} = \{x \in p^{\NN} \st x_i > x_j\}$. Then $\mu(\X_{i,j}) = \frac 1 2 (1- \frac 1 p)$ for all $(i,j)\in \comb {\NN} 2$ and yet for each $x\in \Omega$ the graph $\X(x)=\{(i,j) \in \comb {\NN} 2 \st x_i>x_j\}$ has no paths of length $\geq p$ (where the length of a path is the number of its edges). 
\eexa

 We will next show that the bounds in Example \ref{ex4} are optimal. 
We need:

\bl \label{teoromito}
Let $p\in \NN$ and let $m\in {\mathcal M}^1 (p^{\NN})$. Let 
\begin{equation}
A_{i,j} := \{x \in p^{\NN} \st x_i>x_j \}\,.
\end{equation} 
Then
\begin{equation}\label{mmm}
\inf_{(i,j)\in \comb {\NN} 2} m(A_{i,j}) \leq \frac 1 2 \left( 1 - \frac 1 p \right)
\end{equation}
\el

\bp The proof is a reduction to the case of \emph{exchangeable measures} (see Appendix \ref{app-exch}). Note that if $\sigma \in \Ic$, then $(\sigma \cdot m)(A_{i,j}) = m(A_{\sigma(i), \sigma(j)})$. 
Hence,  replacing $m$ with $\sigma \cdot m$  in \eqref{mmm}  can only increase the infimum, as it is equivalent to the infimum of $m(A_{i,j})$ over a subset of   $ \comb {\NN} 2$. 
By Theorem \ref{lemexch} we can then assume that $m$ is asymptotically exchangeable,
so that in particular the sequence $m_k={\sf S}^k \cdot m$ converges, in the weak$^*$ topology, 
to an exchangeable measure $m^\prime \in {\ca M}^1(p^{\NN})$. 
Since $p$ is finite, the sets $A_{i,j}$ are clopen, and therefore $\lim_{k \to \infty} m_k (A_{i,j}) = m^\prime (A_{i,j}) = m^\prime (A_{0,1})$. Noting that $m_k (A_{i,j}) = m(A_{i+k,j+k})$, it follows that 
\begin{eqnarray}
\inf_{(i,j)\in \comb {\NN} 2} m(A_{i,j}) & \leq & \lim_{k\to \infty} m_k (A_{0,1}) \\
\nonumber & = & m^\prime (A_{0,1}) \\
\nonumber & = & \frac 1 2 \left(1 - m^\prime \{x \st x_0 = x_1\} \right) \\
\nonumber & \leq & \frac 1 2 \left(1 - \frac 1 p \right)
\end{eqnarray}
where the latter inequality follows from Corollary \ref{1p}.
\ep

\bt \label{finpath} Let $(\Omega, \ca A, \mu)$ be a probability space and let $\X : \Omega \to 2^{E_G}$ be a random subgraph of $G:= (\NN, \comb {\NN} 2)$. Consider the set
\[P := \{x \in \Omega \st \X(x) \mbox{ has a path of length } \geq p \}.\] 
Assume $\inf_{e\in \comb {\NN} 2} \mu(\X_e) > \frac 1 2 (1 - \frac 1 p)$. Then 
$\mu(P)>0$. 
\et

\noindent A different proof of this result has been given in \cite[3F]{FT:85} (when the probability space $\Omega$ is $[0,1]$ equipped with the Lebesgue measure).

\bp Suppose for a contradiction that $\mu(P)=0$. We can then assume $P=\emptyset$ (otherwise replace $\Omega$ with $\Omega - P$). For $x\in \Omega$ let $\varphi(x) \colon \NN \to p$ assign to each $i\in \NN$ the length of the longest path starting from $i$ in $\X(x)$. We thus obtain a function $\varphi \colon \Omega \to p^{\NN}$ which is easily seen to be measurable (this is a special case of Lemma \ref{meas}). Let $m = \varphi_\#(\mu) \in {\ca M}^1(p^{\NN})$. Since $\varphi(\X_{i,j}) \subset A_{i,j}$, we have $m(A_{i,j}) \geq \mu(\X_{i,j}) \geq 1/2(1 - \frac 1 p)$ for all $i,j$, contradicting Lemma \ref{teoromito}.  
\ep 

Having determined the critical threshold $\lambda_p= \frac 1 2 (1 - \frac 1 p)$, it follows that if $\inf_{e\in \comb {\NN} 2} \mu(\X_e) \geq \lambda \geq \lambda_p$, the lower bound for $\mu(P)$ grows linearly with $\lambda$. More precisely we have: 

\bc \label{conditional} In the setting of Theorem \ref{finpath}, let $\lambda \in [0,1]$ and suppose that
$\inf_{e\in \comb {\NN} 2} \mu(\X_e) \geq \lambda$. Then $\mu(P) \geq \frac {\lambda - \lambda_p} {1 - \lambda_p}$ where $\lambda_p = \frac 1 2 (1 - \frac 1 p)$. \ec

\bp Suppose $\inf_{e\in \comb {\NN} 2} \mu(\X_e) \geq \lambda$. Consider the conditional probability $\mu(\cdot \mid \Omega - P) \in {\ca M}^1(\Omega)$. We have 
\begin{eqnarray}
\mu(\X_e \mid \Omega - P) & \geq & \frac {\mu(\X_e) - \mu(P)} {1 - \mu(P)} \\
\nonumber & \geq & \frac {\lambda - \mu(P)} {1 - \mu(P)}\,.
\end{eqnarray}
Clearly $\mu(P \mid \Omega - P) = 0$. Applying Theorem \ref{finpath} to $\mu(\cdot \mid \Omega - P)$ it then follows that $\frac {\lambda - \mu(P)} {1 - \mu(P)} \leq \lambda_p$, or equivalently $\mu(P) \geq \frac {\lambda - \lambda_p} {1 - \lambda_p}$. 
\ep


\section{Infinite paths}  \label{secmono}


By Theorem \ref{finpath}, if $\inf_{e\in \comb {\NN} 2} \mu(\X_{i,j}) \geq 1/2$, then the random subgraph $\X$ of $(\NN, \comb {\NN} 2)$ has arbitrarily long finite paths, namely for each $p$ there is $x \in \Omega$ (depending on $p$) such that $\X(x)$ has a path of length $\geq p$. We want to show that for some $x\in \Omega$, $\X(x)$ has an infinite path. To this aim it is not enough to find a single $x$ that works for all $p$. Indeed, $\X(x)$ could have arbitrarily long finite paths without having an infinite path. The existence of infinite paths can be neatly expressed in terms of the following definition. 

\bd Let $G$ be a countable directed graph and let $\omega_1$ be the first uncountable ordinal. We recall that the {\em rank function} $\phi_G\colon V_G \to \omega_1 \cup \{\infty\}$ of $G$ is defined as follows. For $i\in V_G$, 
\[\phi_G(i) = \sup_{j: (i,j)\in E_G} \big(\phi_G(j)+1\big).\] This is a well defined countable ordinal if $G$ has no infinite paths starting at $i$. In the opposite case we set \[\phi_G(i) = \infty\] where $\infty$ is a conventional value bigger than all the countable ordinals. For notational convenience we will take $\infty = \omega_1$ so that $\omega_1 \cup \{\infty\}=\omega_1 \cup \{\omega_1\} = \omega_1+1$.
Note that if $i$ is a leaf, $\phi_G(i) = 0$. 
Also note that $G$ has an infinite path if and only if $\phi_G$ assumes the value $\infty$. 

Given a random subgraph $\X \colon \Omega \to 2^{E_G}$ of $G$, we let $\phi_{\X}(x)= \phi_{\X (x)}$, namely $\phi_{\X} (x)(i)$ is the rank of the vertex $i$ in the graph $\X(x)$. So $\phi_{\X}$ is a map from $\Omega$ to $(\omega_1+1)^{V_G}$. It can also be considered as a map from $\Omega \times V_G$ to $\omega_1 + 1$ by writing $\phi_{\X}(x,i)$ instead of $\phi_{\X} (x)(i)$. \ed 

\br \label{ind} 
We have $\phi_{\X}(x,i) = \phi_{\omega_1}(x,i)$ where $\phi_\alpha \colon \Omega \to (\omega_1+1)^{V_G}$ is   
the truncation $\phi_\alpha:=\min (\phi, \alpha)$, that we can equivalently define by induction on $\alpha \leq \omega_1$ as follows. 
\[
\begin{array}{lll} \label{defrec} 
\phi_0(x,i) &=& 0 \\
\phi_\alpha(x,i) &=& \sup \{ \phi_{\beta}(x,j)+1 \st \beta<\alpha,\ (i,j)\in \X(x) \}
\end{array}
\]
\er
The above representation will be of use in the following lemma in connection to measurability properties of the map $\phi$.

\bl \label{meas} Let $G$ be a countable directed graph, let $(\Omega, \ca A, \mu)$ be a probability space and let $\X : \Omega \to 2^{E_G}$ be a random subgraph of $G$. 
\begin{enumerate}
\item The set $P:= \{x \in \Omega \st \X(x) \text{ has an infinite path } \}$ is $\mu$-measurable. 
\item For all $\alpha \leq \omega_1$ and $i\in V_G$, the set $\{x \in \Omega \st \phi_{\X}(x,i) = \alpha\}$ is $\mu$-measurable. 
\item $\varphi_{\X}\colon \Omega \to {(\omega_1+1)}^{V_G}$ is $\mu$-measurable and its  restriction to $\Omega-P$ is essentially bounded, namely for some $\alpha_0<\omega_1$ it takes values in $\alpha_0^{V_G}$ outside of a $\mu$-null set. 
\end{enumerate} 
\el 
\bp 
Since taking the supremum over a countable set preserves measurability, 
from Remark \ref{ind} it follows that for all $i\in V_G$ and $\alpha < \omega_1$ the sets 
$\{x: \phi_{\X}(x,i)=\alpha\}$ are measurable. We will show that $\{x: \phi_{\X}(x,i)=\omega_1\}$ is $\mu$-measurable, namely it is the union of a measurable set and a $\mu$-null set. Fix $i\in V_G$. The sequence of values
$\mu\left(\left\{ x:\,\phi_{\X}(x,i)\le\beta\right\}\right)$ is increasing with respect to the countable ordinal $\beta$ and uniformly bounded by $1=\mu(\Omega),$ therefore it is stationary at some finite value. So there is $\alpha_0<\omega_1$ such that 
\[
\mu\left(\left\{ x\in \Omega:\ \phi_{\X}(x,k)=\beta\right\}\right)=0 
\qquad \mbox{ for } \alpha_0\le\beta<\omega_1\,. 
\]
It follows that $\{x: \phi_{\X}(x,i)=\omega_1\}$ is $\mu$-measurable and $\varphi_{\X}$ is $\mu$-measurable. Since $P = \cup_i \{x: \phi_{\X}(x,i) = \omega_1\}$, we have that $P$ is $\mu$-measurable,too. 
\ep

Given an ordinal $\alpha$, we put on $\alpha$ the topology generated by the open intervals. Note that a non-zero ordinal is compact if and only if it is a successor ordinal, and it is metrizable if and only if it is countable. Let ${\mathcal M}_c (\omega_1^{\NN})$ be the set of compactly supported Borel measures on $\omega_1^{\NN}$, namely the measures with support in ${\alpha_0}^{\NN}$ for some ${\alpha_0} < \omega_1$. The following Lemma reduces to Lemma \ref{teoromito} if $\alpha_0$ is finite. 

\bl \label{teoromito2}
Let $m\in {\mathcal M}_c (\omega_1^{\NN})$ be a non-zero measure with compact support. Let 
\begin{equation}
A_{i,j} := \{x \in \omega_1^{\NN} \st x_i>x_j \}\,.
\end{equation} 
Then
\begin{equation}\label{mmm2}
\inf_{(i,j)\in \comb {\NN} 2}m(A_{i,j})<\frac {m\left(\omega_1^{\NN}\right)}{2}\,.
\end{equation}
\el

\bp 
With no loss of generality we can assume that $m\in \mathcal M^1(\omega_1^{\NN})$,
i.e. $m\left(\omega_1^{\NN}\right)=1$. 
We divide the proof into four steps.

\noindent 
{\it Step 1.} Letting $\partial\omega_1$ be the derived set of $\omega_1$, that is the subset 
of all countable limit ordinals, we can assume that
\[
m\left( \left\{ x:\ x_i\in\partial\omega_1\right\}\right) = 0
\qquad \forall i\in\NN .
\]
Indeed, it is enough to observe that the left-hand side of equation \eqref{mmm2} can only increase if
we replace $m$ with $s_\#(m)$, where $s:\,\omega_1\to \omega_1\setminus \partial\omega_1$
is the successor map sending $\alpha<\omega_1$ to $\alpha+1$, and $s_\#(m) = (s_*)_\#$, namely $s_\#(m)(X):= m(\{ x \in \omega_1^{\NN} \st s\circ x \in X\})$.  

\noindent 
{\it Step 2.} Since the support of $m$ is contained in $\alpha_0^{\NN}$, for some ordinal $\alpha_0<\omega_1$,
thanks to Theorem \ref{lemexch} we can assume that $m$ is asymptotically exchangeable,
i.e. the sequence $m_k={\sf S}^k \cdot \sigma\cdot m$ converges, in the weak$^*$ topology, 
to an exchangeable measure $m^\prime \in {\ca M}^1(\omega_1^{\NN})$, with support in $\alpha_0^{\NN}$, for all $\sigma\in\omega^{(\omega)}$. 
Note however that, unless $\alpha_0$ is finite, we cannot conclude  that $\lim_{k\to \infty} m_k(A_{i,j}) = m'(A_{i,j})$ since the sets 
$A_{i,j} = \{x \in \omega_1^{\NN} \st x_i > x_j \}$ are not clopen. 

\noindent 
{\it Step 3.} We shall prove by induction on $\alpha<\omega_1$ that
\begin{equation}\label{ggg}
\liminf_{(i,j)\to +\infty} m\left( \left\{ x:\ x_j<x_i\le\alpha\right\}\right)\le 
m^\prime\left( \left\{ x:\ x_1<x_0\le\alpha\right\}\right).
\end{equation}
Indeed, for $\alpha=0$ we have $\left\{ x:\ x_j<x_i\le 0\right\}=\emptyset$, 
and \eqref{ggg} holds.

As inductive step, let us assume that \eqref{ggg} holds for all $\alpha<\beta<\omega_1$,
and we distinguish whether $\beta$ is a successor or a limit ordinal. 

In the former case let $\beta=\alpha +1$.  
For $(i,j) \to +\infty$ (with $i<j$) we have:
\begin{eqnarray*}
m \left(\left\{ x_j< x_i\le \beta\right\}\right) &=&
m \left(\left\{ x_j< x_i \le \alpha\right\}\right) +
m \left(\left\{ x_j\le\alpha,\  x_i= \beta\right\}\right)
\\
& \leq & 
m^\prime\left(\left\{ x_1<x_0\le \alpha\right\}\right) \! +
m^\prime\left(\left\{ x_1\le\alpha,\  x_0= \beta\right\}\right) \! +o(1)
\\
&=& m^\prime\left(\left\{ x_1<x_0\le \beta\right\}\right) + o(1)\,,
\end{eqnarray*}
where we used the induction hypothesis, and the fact that 
$\left\{ x_j\le\alpha,\  x_i= \beta\right\}$ is clopen.

Let us now assume that $\beta$ is a limit ordinal and let $i \in \NN$. We have
\[
\bigcap_{\alpha<\beta}\left\{ x:\ \alpha<x_i<\beta\right\}=\emptyset ,
\]
so for all $\varepsilon > 0$ there exists $\alpha<\beta$ such that 
\[
m^\prime\left( \left\{ \alpha<x_i<\beta\right\}\right)<\eps \,.
\]
Since $m^\prime$ is exchangeable, we can choose the same $\alpha$ for every $i$. Moreover by assumption $m(\{x_i = \beta\}) = 0$ for every $i\in \NN$. Hence there exists $\alpha \leq \alpha_i<\beta$ such that 
\[
m(\{\alpha_i \leq x_i \leq \beta\}) < \eps\,.
\]
Given $i<j$, distinguishing the relative positions of $x_i,x_j$ with respect to $\alpha$ and $\alpha_i$ we have:  
\begin{eqnarray*}
\left\{ x_j< x_i\le \beta\right\} &\subseteq&
\left\{ x_j< x_i\le \alpha\right\} 
\\ && \cup 
\left\{ x_j\le\alpha < x_i\le \beta\right\}
\\ && \cup 
\left\{ \alpha< x_j\le \alpha_i\right\}
\\ && \cup 
\left\{ \alpha_i< x_i\le \beta\right\}\,.
\end{eqnarray*}
which gives
\begin{eqnarray} \label{sums} 
m\left(\left\{ x_j< x_i\le \beta\right\}\right) &\le&
m\left(\left\{ x_j< x_i\le \alpha\right\}\right) 
\\ \nonumber && +
m\left(\left\{ x_j\le\alpha < x_i\le \beta\right\}\right)
\\ \nonumber && +
m\left(\left\{ \alpha< x_j\le \alpha_i\right\}\right)
\\ \nonumber && +
m\left(\left\{ \alpha_i< x_i\le \beta\right\}\right).
\end{eqnarray}

Since $\left\{ x_j\le\alpha < x_i\le \beta\right\}$ and $\left\{ \alpha< x_j\le \alpha_i\right\}$ are both clopen, we can approximate their $m$-measure by their $m^\prime$-measure. So we have:

\begin{eqnarray*}
m\left\{ x_j\le\alpha < x_i\le \beta\right\} & = & m^\prime \left(\left\{ x_1\le\alpha < x_0\le \beta\right\}\right) + o(1) \\
& &  \mbox{ for } (i,j) \to \infty 
\end{eqnarray*}
and
\begin{eqnarray*}
m\left(\left\{ \alpha< x_j\le \alpha_i\right\}\right) & =&  m^\prime \left(\left\{ \alpha< x_1 \le \alpha_i\right\}\right) +o(1) \\
 & & \mbox{ for } j\to \infty,
\end{eqnarray*}
where we used Remark \ref{j} to allow $j\to \infty$ keeping $i$ fixed. 
Now note that by the choice of $\alpha$, we have $m^\prime \left(\left\{ \alpha< x_1 \le \alpha_i\right\}\right) <\eps$, and by induction hypothesis $\liminf_{(i,j)\to +\infty} m \left(\left\{ x_j< x_i\le \alpha\right\}\right) < m'(\{x_1<x_0 \leq \beta\})$. Hence, from \eqref{sums} we obtain: 
\begin{eqnarray*}
\liminf_{(i,j)\to +\infty} m\left(\left\{ x_j< x_i\le \beta\right\}\right) 	&\le&
m'(\{x_1<x_0 \leq \alpha \})
\\ && +
m^\prime \left(\left\{ x_1\le\alpha < x_0\le \beta\right\}\right) 
\\ && +
\eps + \eps \,.
\end{eqnarray*}
Therefore, 
\begin{eqnarray*}
\liminf_{(i,j)\to +\infty} 
m\left(\left\{ x_j< x_i\le \beta\right\}\right)  &\le & 
m^\prime\left(\left\{ x_1< x_0\le \beta\right\}\right) + 2\eps 
\end{eqnarray*}
Inequality \eqref{ggg} is then proved for all $\alpha<\omega_1$.

\noindent 
{\it Step 4.} We now conclude the proof of the theorem. 
From \eqref{ggg} it follows
\begin{equation}\label{labella}
\inf_{(i,j)\in \comb {\NN} 2}m\left( A_{i,j}\right)\le 
m^\prime\left( \left\{ x:\ x_1<x_0\right\}\right) = 
\frac 1 2 \left( 1- m^\prime\left( \left\{ x:\ x_1=x_0\right\}\right)\right)
< \frac 1 2 \,.
\end{equation}
where we used the fact the $m'$ is exchangeable and Corollary \ref{zero}. 
\ep

\bt \label{infpath} Let $(\Omega, \ca A, \mu)$ be a probability space and let $\X : \Omega \to 2^{E_G}$ be a random subgraph of $G:= (\NN, \comb {\NN} 2)$. Consider the set
\[
P := \{x \in \Omega \st \X(x)\ \mbox{\rm has an infinite path} \}.
\] 
Assume $\inf_{e\in \comb {\NN} 2} \mu(\X_e) \geq \frac 1 2$. Then $\mu(P)>0$.  
\et 

\noindent As observed in the Introduction, we recall that this result  follows from \cite[4D]{FT:85}, when $\Omega=[0,1]$ 
with the Lebesgue measure.

\bp
Suppose for a contradiction $\mu(P)=0$. We can then assume $P=\emptyset$ (replacing $\Omega$ with $\Omega - P$). Hence the rank function $\varphi := \varphi_{\X}\colon \Omega \to (\omega_1 + 1)^{\NN}$ takes values in $\omega_1^{\NN}$. Let $m = \varphi_\#(\mu) \in {\ca M}^1(\omega_1^{\NN})$. Note that $\varphi(\X_{i,j}) \subset A_{i.j} := \{x \in p^{\NN} \st x_i> x_j \}$. Hence $m(A_{i,j}) \geq \mu(\X_{i,j}) \geq 1/2$ for all $(i,j)\in \comb {\NN} 2$. This  contradicts Lemma \ref{teoromito2}.  
\ep 

\br
Note that the bound $1/2$ is optimal by Example \ref{ex4}. 
\er

Reasoning as in Corollary \ref{conditional} we obtain: 

\bc \label{conditional2} 
Let $0\leq \lambda < 1$. If $\inf_{e\in \comb {\NN} 2} \mu(\X_e) \geq \lambda$, then $\mu(P) > \frac {\lambda - 1/2} {1 - 1/2}$. \ec

Note that if we replace $(\NN, \comb {\NN} 2)$ with a finitely branching countable graph $G$, then the threshold for the existence of infinite paths becomes $1$, namely we cannot ensure the existence of infinite paths even if each edge of $G$ belongs to the random subgraph $\X$ with probability very close to $1$. In fact, the following more general result holds: 

\bprop \label{fin-b} Let $G = (V_G, E_G)$ be graph admitting a coloring function $c \colon E_G\to \NN$ such that each infinite path in $G$ meets all but finitely many colours (it is easy to see that a finitely branching countable graph $G$ has this property). Then for every $\eps > 0$ there is a probability space $(\Omega, \ca A, \mu)$ and a random subgraph $\X \colon \Omega \to 2^{E_G}$ of $G$ such that for all $x\in \Omega$, $\X(x)$ has no infinite paths, and yet $\mu(\X_e)> 1-\eps$ for all $e\in E_G$. 
\eprop

\bp Let $(Z_n)_{n\in \NN}$ be a disjoint family of infinite subsets of $\NN$. Let $\mu$ be a probability measure on $\Omega:=\NN$ with $\mu(\{n\}) < \eps$ for every $n$. Given $n\in \Omega$ let $\X(n)$ be the subgraph of $G$ (with vertices $V_G$) containing all edges $e\in E_G$ of colour $c(e) \nin Z_n$. Given $e\in E_G$ there is at most one $n$ such that $c(e)\in Z_n$. Hence clearly $\mu(\X_e) \geq 1-\eps$, and yet $\X(n)$ has no infinite paths for any $n\in \Omega$.
\ep

\br \label{reals} 
It is natural to ask whether the answer to Problem \ref{increasing} changes 
if we substitute $\NN$ with the set of the real numbers. Since $\NN\subset\R$, the 
probability threshold for the existence of infinite paths can only decrease, 
but the following example shows that it still equals $1/2$.
Let $\Omega=[0,1]^{\R}$ equipped with the product Lebesgue measure $\mathcal L$, let $\eps>0$, and let 
\[
\X_{i,j} := \big\{ x\in\Omega:\, x_i>x_j+\eps\big\}\,, 
\]
for all $i<j\in\R$.
The assertion follows observing that $\mathcal L(\X_{i,j})= (1-\eps)^2/2$ for all $i<j\in\R$, 
and 
\[
\bigcap_{i\in\{1,\ldots,N\}}\X_{n_i,n_{i+1}}=\emptyset
\]
whenever $n_i$ is a strictly increasing sequence of real numbers, and $N>1/\eps$.  
\er

\section{Threshold functions for graph morphisms} \label{secrelp}

\bd Let $F$ and $G$ be directed graphs. A graph morphism $\p \colon G \to F$ is a map $\p:\, V_G\to V_F$ such that $(\p(a),\p(b))\in E_F$ for all $(a,b)\in E_G$. We write $G\to F$ if there is a graph morphism from $G$ to $F$. \ed 

The results of the previous sections were implicitly based on following observation:

\br \label{implicit} Let $G$ be a directed graph. 
\begin{enumerate}
\item 
$G$ has a path of length  $\geq p$ if and only if $G \not\to (p, \comb p 2)$.
\item $G$ has an infinite path if and only if $G\not\to (\omega_1, \comb {\omega_1} 2)$.  
\end{enumerate} 
\er

This suggests to generalize the above results considering other properties of graphs that can be expressed in terms of non-existence of graph morphisms. 
Let us give the relevant definitions.

\bd \label{relcap} 
Given two directed graphs $F,G$ and given $i,j\in V_G$ let 
\begin{equation}
A_{i,j}(F,G) := \{ u\in V_F^{V_G} :\, (u(i),u(j))\in E_F\}
\end{equation} 
and define the \emph{relative capacity} 
of $F$ with respect to $G$ as
\begin{equation}\label{eqcap}
c(F,G) := \sup_{m\in {\mathcal M}^1(V_F^{V_G})}\; \inf_{(i,j)\in E_G}
m\left(A_{i,j}\left(F,G\right)\right) \in [0,1].
\end{equation}
\ed

%

Theorems \ref{finpath} and \ref{infpath} have the following counterpart. 

\bt \label{capacity} Let $F$ and $G$ be directed countable graphs, let $(\Omega, \ca A, \mu)$ be a probability space and let $\X\colon \Omega \to 2^{E_G}$ be a random subgraph of $G$. Let $P := \{x\in \Omega \st \X(x) \not\to F \}$. 
Assume $\inf_{e\in E_G} \mu(\X_e) > c(F,G)$. Then $\mu(P)>0$. 
Moreover there are examples in which $P$ is empty and $\inf_{e\in E_G} \mu(\X_e)$ is as close to $c(F,G)$ as required. So $c(F,G)$ is the threshold for non-existence of graph morphisms $f\colon \X(x)\to F$. To prove the second part it suffices to take $\Omega = V_F^{V^G}$ and $\X_{i,j} = A_{i,j}(F,G)$. 
\et

\bp Suppose for a contradiction $\mu(P)=0$. We can then assume $P=\emptyset$ (replacing $\Omega$ with $\Omega - P$). Hence for each $x\in \Omega$ there is a graph morphism $\varphi(x) \colon \X(x)\to F$, which can be seen as an element of $V_F^{V_G}$. We thus obtain a map $\varphi \colon \Omega \to V_F^{V_G}$. By Lemma \ref{measurable} below, $\varphi$ can be chosen to be $\mu$-measurable. Since $x \in \X_{i,j}$ implies $(\phi(x)(i),\phi(x)(j))\in E_F$, we have $\varphi(\X_{i,j}) \subset A_{i,j}(F,G)$ for all $(i,j)\in E_G$. Let $m:= 
\varphi_\#(\mu) \in {\ca M}^1(V_F^{V_G})$.  Then $m(A_{i,j}(F,G)) \geq \mu(\X_{i,j}) > c(F,G)$. This is absurd by definition of $c(F,G)$. 
\ep

Reasoning as in Corollary \ref{conditional} we obtain: 

\bc \label{conditional3} Suppose $c(F,G) <1$. 
If $\inf_{e\in \comb {\NN} 2} \mu(\X_e) \geq \lambda$, then $\mu(P) \geq \frac {\lambda - c(F,G)} {1 - c(F,G)}$. \ec

\br 
If the sup in the definition of $c(F,G)$ is not reached, it suffices to have
the weak inequality $\inf_{e\in E_G} \mu(\X_e) \geq c(F,G)$ in order to have $\mu(P)>0$ (this is indeed the case of Theorem \ref{infpath}). 
\er

It remains to show that the map $\varphi \colon \Omega \to V_F^{V_G}$ in the proof of Theorem \ref{capacity} can be taken to be $\mu$-measurable. 

\bl \label{measurable} Let $F,G$ be countable directed graphs, let $(\Omega, \ca A, \mu)$ be a probability space, and let $\X \colon \Omega \to 2^{E_G}$ be a random subgraph of $G$. 
\begin{enumerate}
\item The set $\Omega_0 := \{x \in \Omega \st \X(x) \to F\}$ is $\mu$-measurable (i.e. measurable with respect to the $\mu$-completion of $\ca A$). 
\item There is an $\mu$-measurable function $\varphi \colon \Omega_0 \to V_F^{V_G}$ that selects, for each $x\in \Omega_0$, a graph morphism $\phi(x)\colon \X(x) \to F$. 
\item If $F$ is finite, then $\Omega_0$ is measurable and $\varphi$ can be chosen measurable.
\end{enumerate}
\el 

\begin{proof} Given a function $f\colon V_G \to V_F$, we have $f\colon \X(x) \to F$ (i.e., $f$ is a graph morphism from $\X(x)$ to $F$) if and only if $x \in \bigcap_{(i,j)\in V_G}\bigcup_{(a,b)\in V_F} B_{i,j,a,b}$, where $x\in B_{i,j,a,b}$ says that $f(i)=a, f(j)=b$ and $x\in \X_{i,j}$. This shows that $B := \{(x,f) \st f \colon \X(x)\to F\}$ is a measurable subset of $\Omega \times V_G^{V_F}$. We are looking for a ($\mu$-)measurable function $\varphi \colon \pi_X(B)\to V_F^{V_G}$ whose graph is contained in $B$. 

Special case: Let us first assume that $\Omega$ is a Polish space (i.e., a complete separable metric space) with its algebra $\ca A$ of Borel sets. By Jankov - von Neumann uniformization theorem (see \cite[Thm. 29.9]{K:95}), if $X,Y$ are Polish spaces and $Q\subset X\times Y$ is a Borel set, then the projection $\pi_X(Q) \subset X$ is universally measurable (i.e. it is $m$-measurable for every $\sigma$-finite Borel measure $m$ on $X$), and there is a universally measurable function $f\colon \pi_X(Q) \to Y$ whose graph is contained in $Q$. We can apply this to $X=\Omega, Y=V_F^{V_G}$ and $Q=B$ to obtain (1) and (2). It remains to show that if $F$ is finite $\pi_X(Q)$ and $f$ can be chosen to be Borel measurable. To this aim it suffices to use the following uniformization theorem of Arsenin - Kunugui (see \cite[Thm. 35.46]{K:95}): if $X,Y,Q$ are as above and each section $Q_x = \{y\in Y\st (x,y)\in Q\}$ is a countable unions of compact sets, then $p_X(Q)$ is Borel and there is a Borel measurable function $f\colon \pi_X(Q) \to Y$ whose graph is contained in $Q$. 

General case: We reduce to the special case as follows. Let $X = 2^{V_G}, Y = V_F^{V^G}$ and
consider the set $B' \subset X \times Y$ consisting of those pairs $(H,f)$ such that $H$ is a subgraph of $G$ (with the same vertices) and $f\colon H\to F$ is a graph morphism. Consider the pushforward measure $m = \X_\# (\mu)$ defined on the Borel algebra of $2^{V_G}$. By the special case there is a ($m$-)measurable function $\psi \colon \pi_X(B')\to V_F^{V_G}$ whose graph is contained in $B'$. To conclude it suffices to take $\varphi := \psi \circ \X$.
\end{proof} 

We now show how to compute the relative capacity $c(F, (\NN, \comb {\NN} 2))$ (see Definition \ref{relcap}) for any finite graph $F$. 
The following invariant of directed graphs has been studied in \cite{Rao:82} and  
\cite[Section 3]{FT:85}.

\bd 
Given a directed graph $F$, we define the \emph{capacity} of $F$ as
\begin{equation}\label{eqcoop}
c_0(F) := \sup_{\lambda\in \Sigma_F} \; \sum_{(a,b)\in E_F}\lambda_a\lambda_b
\quad \in \,[0,1]\,,
\end{equation}
where $\Sigma_{F}$ is 
the symplex of all sequences $\{\lambda_a\}_{a\in V_F}$ of real numbers such that 
$\lambda_a\ge 0$ and $\sum_{a\in V_F}\lambda_a=1$.
\ed 

\bprop\label{puno}
If $F$ is a finite directed graph, then
\begin{equation}\label{comcap}
c\left(F, (\NN,\comb {\NN} 2)\right) = c_0(F) \,.
\end{equation}
\eprop
\bp Let $G = (\NN, \comb {\NN} 2)$. The proof is a series of reductions. 

\noindent {\it Step 1.}
Note that if $\sigma \in \Ic$, then $\sigma \cdot m(A_{i,j}(F,G)) = m(A_{\sigma(i), \sigma(j)})$. Hence the infimum in \eqref{eqcap} can only increase replacing $m$ with $\sigma^*_\#(m)$. By Theorem \ref{lemexch} there is $\sigma \in \Ic$ such that $\sigma \cdot m$ is asymptotically exchangeable. It then follows that 
we can equivalently take the supremum in \eqref{eqcap} among the measures 
$m\in {\mathcal M}^1(V_F^{\NN})$ which are asymptotically exchangeable.  

\noindent {\it Step 2.} By definition if $m$ is asymptotically exchangeable there is an exchangeable measure $m'$ such that $\lim_{k\to \infty} m_k = m'$, where $m_k = {\sf S}^k\cdot m$. Clearly $$\inf_{(i,j)\in E_G } m (A_{i,j}(F,G)) \leq \lim_{k\to \infty} m_k (A_{0,1}(F,G)) = m'(A_{0,1}(F,G)).$$ So the supremum in \eqref{eqcap} coincides with $\sup_m m\left(A_{0,1}\left(F,G\right)\right)$, 
for $m$ ranging over the exchangeable  measures. 

\noindent {\it Step 3.}
Recalling \eqref{choquetbern}, 
every exchangeable measure is a convex integral combination of Bernoulli 
measures $B_\lambda$, with $\lambda\in \Sigma_F$. It follows that 
it is sufficient to compute the supremum on the 
Bernoulli measures $B_\lambda$. We have: 
\begin{eqnarray*} 
B_\lambda \left(\left\{ x \in V_F^{\NN} \st (x_0, x_1) \in E_F \right\} \right) & = & \sum_{(a,b)\in E_F} B_\lambda \left(\left\{x \st x_0 = a \;, \; x_1 = b \right\} \right) \\
& = & \sum_{(a,b) \in E_F} \lambda_a \lambda_b
\end{eqnarray*}
so that \eqref{eqcap} reduces to \eqref{eqcoop}.
\ep

Notice that if there is a morphism of graphs from $G$ to $F$, then 
$c_0(G) \leq c_0(F)$. Also note that $c_0(F)=1$ if there is some $a\in V_F$ with 
$(a,a)\in E_F$. Recall that $F$ is said to be: {\em irreflexive} if $(a,a)\nin E_F$ for all $a\in V_F$; {\em symmetric} if $(a,b)\in E_F \Longleftrightarrow (b,a) \in E_F$ for all $a,b\in V_F$; {\em anti-symmetric} if $(a,b)\in E_F \Longrightarrow (b,a)\nin E_F$ for all $a,b\in V_F$. 

The {\em clique number} ${\rm cl}(F)$ of $F$ is defined as the largest integer $n$ such that there is a subset $S \subset V_F$ of size $n$ which forms a clique, namely $(a,b)\in E_F$ or $(b,a)\in E_F$ for all $a,b\in S$. 

\bprop\label{ppie} {\rm(see also \cite[Section 3]{FT:85})}
Let $F$ be a finite irreflexive directed graph. If $F$ is anti-symmetric, then 
\begin{equation}\label{cappa}
c_0(F) =  \frac 1 2 \left( 1-\frac{1}{{\rm cl}(F)} \right) \,.
\end{equation}
If $F$ is symmetric, then 
\begin{equation}
c_0(F) =  1-\frac{1}{{\rm cl}(F)}\,.
\end{equation}
In particular $c_0(K_p) = 1-\frac{1} p$. 
\eprop

\bp The anti-symmetric case follows from the symmetric one taking the symmetric closure. So we can assume that $F$ is symmetric. 
Let $\lambda\in\Sigma_F$ be a maximizing distribution, 
meaning that $c_0(F) = \sum_{(a,b)\in E_F}\lambda_a\lambda_b$,
and let $S_\lambda$ be the subgraph of $F$ spanned by the support of $\lambda$, 
that is $V_{S_\lambda}=\{ a\in V_F:\, \lambda_a>0\}$. Given $a\in S_\lambda$ note that $\frac \partial {\partial \lambda_a} \sum_{(u,v)\in E_F} \lambda_u \lambda_v = 2 \sum_{b \in V_F:\, (a,b)\in E_F } \lambda_b$. 
{}From Lagrange's multiplier Theorem it then follows that $\sum_{b \in V_F:\, (a,b)\in E_F } \lambda_b$ is constant, namely it does not depend on the choice of $a\in S_\lambda$. Since $\sum_{a\in S_\lambda} (\sum_{b\st (a,b) \in E_F} \lambda_a) = c_0(F)$, it follows that for each $a\in S_\lambda$ we have:
\begin{equation} \label{eqlb}
\sum_{b \in V_F:\, (a,b)\in E_F } \lambda_b = c_0(F)\,.
\end{equation}
If $c,c^\prime\in V_{S_\lambda}$, we can consider the distribution
$\lambda^\prime\in \Sigma_F$ such that 
$\lambda^\prime_c=0$, $\lambda^\prime_{c^\prime}=\lambda_c+\lambda_{c^\prime}$, and 
$\lambda^\prime_b=\lambda_b$ for all $b\in V_F\setminus\{c, c^\prime\}$.
{}From \eqref{eqlb} it then follows that $\lambda^\prime$ is also a maximizing 
distribution whenever $(c,c')\nin E_F$. (In fact $\sum_{(a,b)\in E_F} \lambda'_a \lambda'_b = \sum_{(a,b)\in E_F} \lambda_a \lambda_b - \lambda_c \sum_{b \st (c,b)\in E_F} \lambda_b +  \lambda _{c} \sum_{b \st (c',b)\in E_F} \lambda b = c_0(F) - \lambda_c c_0(F) + \lambda_c c_0(F)$.)  

As a first consequence, $S_\lambda$ is a clique whenever $\lambda$ is a maximizing distribution with minimal support.
Indeed, let $K$ be a maximal clique contained in $S_\lambda$, and assume by 
contradiction that there exists $a\in V_{S_\lambda}\setminus V_K$. Letting $a^\prime\in V_K$ 
be a vertex of $F$ independent of $a$  
(such an element exists since $K$ is a maximal clique), and letting $\lambda^\prime\in \Sigma_F$ as above,  
we have $c_0(F) = \sum_{(a,b)\in E_F}\lambda^\prime_a\lambda^\prime_b$, contradicting the minimality of $V_{S_\lambda}$.

Once we know that $S_\lambda$ is a clique, 
again from \eqref{eqlb} we get that $\lambda$ is a uniform ditribution, that is 
$\lambda_a=\lambda_b$, for all $a,b\in V_{S_\lambda}$. It follows
\[
c_0(F) \; = \; 1- \frac{1}{\vert S_\lambda\vert} 
\; \le \; 1-\frac{1}{{\rm cl}(F)} \,,
\]
which in turn implies \eqref{cappa}, the opposite inequality being realized by 
a uniform distribution on a maximal clique.
\ep

Notice that the proof of Proposition \ref{ppie} shows that there exists a 
maximizing $\lambda\in\Sigma_F$ whose support is a clique 
(not necessarily of maximal order).

\subsection{Chromatic number} 

We will apply the results of the previous section to study the chromatic number of a random subgraph of $(\NN, \comb {\NN} 2)$. 
We point out that an alternative proof of this result follows from \cite[Theorem 1]{EH:64}.

We recall that the chromatic number $\chi(G)$ of a directed graph $G$ is the smallest $n$ such that there is a colouring of the vertices of $G$ with $n$ colours in such a way that $a,b\in V_G$ have different colours whenever $(a,b)\in E_G$ (see \cite{Bollobas}). 

For $p\in \NN$, let $K_p$ be the complete graph on $p$ vertices, namely $K_p$ has set of vertices $p=\{0,1,\ldots, p-1\}$ and set of edges $\{(x,y)\in p^2 \st x\neq y\}$. 
Clearly $\chi(K_p) = p$. Note also that: 

\begin{equation}\label{cr}
G\to K_p \Longleftrightarrow \chi(G)\le p\,.
\end{equation}

Now let $(\Omega, \ca A, m)$ be a probability space, and let $\X \colon \Omega \to 2^{E_G}$ be a random subgraph of $G = (\NN, \comb {\NN} 2)$. Let $P = \{x\in \Omega \st \chi(\X(x)) \geq p\}$. 
By Equation \eqref{cr} and the results of the previous section, if $\inf_{e\in \mu(\X_e)} > c(K_p, (\NN, \comb {\NN} 2)$, then $\mu(P)>0$. This however does not say much unless we manage to determine $c(K_p, (\NN, \comb {\NN} 2))$. We will show that $c(K_p, (\NN, \comb {\NN} 2)) = (1 - \frac 1 p)$, so we have: 

\bt
Let $(\Omega, \ca A, m)$ be a probability space, and let $\X \colon \Omega \to 2^{E_G}$ be a random subgraph of $(\NN, \comb {\NN} 2)$. If $\inf_{e\in \mu(\X_e)} > 1 - \frac 1 p$, then 
$$
\mu(\{x\in \Omega \st \chi(\X(x)) \geq p+1\}) > 0.
$$ 
\et

\section{Infinite cliques} \label{secproblem}

We recall the following standard Borel-Cantelli type result, which shows that Problem \ref{problemstrong} has a positive answer for $k=1$. 

\bprop \label{lemone} Let $(\Omega, \ca A, \mu)$ be a probability space. 
Let $\lambda > 0$ and for each $i\in \NN$ 
let $X_i\subseteq \Omega$ be a measurable set such that $\mu(X_i)\ge\lambda$.
Then there is an infinite set $J\subset \NN$ such that 
\[
\bigcap_{i\in J}X_i \ne\emptyset.
\] 
\eprop

\bp 
The set $Y:= \bigcap_n \bigcup_{i>n} X_i$ is a decreasing intersection of sets of (finite) measure greater than $\lambda>0$, 
hence $\mu(Y)\ge\lambda$ and, in particular, $Y$ is non-empty. 
Now it suffices to note that any element $x$ of $Y$ belongs to infinitely many $X_i$'s. 
\ep

Proposition \ref{lemone} has the following interpretation: if we choose each element of $\NN$ with probability greater or equal to $\lambda$, we obtain an infinite subset with probability greater or equal to $\lambda$.

The following example shows that Problem \ref{problemstrong} has in general a negative answer for $k>1$. 

\bexa\label{exo}   
Let $p\in \NN$ and consider the Cantor space $\Omega = p^{\NN}$, equipped  
with the Bernoulli measure $B_{(1/p,\ldots,1/p)}$, and let $\X_{i,j} := \{x \in \Omega :\, x_i \neq x_j\}$. 
Then each  $\X_{i,j}$ has measure $\lambda = 1-1/p$, and
for all $x\in X$ the graph 
$\X(x):= \{ (i,j)\in \comb {\NN} 2: x\in \X_{i,j}\}$ does not contains cliques 
(i.e. complete subgraphs) of cardinality $(p+1)$. 
\eexa

In view of Example \ref{exo}, we need further assumptions in order to get a positive answer to Problem \ref{problemstrong}. 

\bexa \label{excountable} By Ramsey theorem, Problem \ref{problemstrong} has a positive answer 
if there is a finite set $S\subset \Omega$ such that each $X_{i_1, \ldots, i_k}$ has a non-empty intersection with $S$. In particular, this is the case if $\Omega$ is countable. 
\eexa

\bprop 
Let $r>0$. Assume that $\Omega$ is a compact metric space and 
each set $\X_{i_1, \ldots, i_k}$ contains a ball $B_{i_1, \ldots, i_k}$ of
radius $r>0$. Then Problem \ref{problemstrong} has a positive answer. 
\eprop

\bp 
Applying Lemma \ref{lempie} to the centers of the balls $B_{i_1, \ldots, i_k}$ 
it follows that for all $0<r'<r$ there exists an infinite set $J$ and a ball $B$ of radius $r'$ such that 
\[
B\subset\bigcap_{(j_1, \ldots ,j_k)\in J^{[k]}}X_{j_1, \ldots, j_k}.
\]  
\ep

We now give a sufficient condition for a positive answer to Problem \ref{problemstrong}.

\bt\label{taiut} Let $(\Omega, \ca, \mu)$ be a probability space. Let $\lambda > 0$ and 
assume that we have the sets $\mu(\X_{i_1 \ldots i_k})\geq \lambda$ for each $(i_1, \ldots, i_k)\in \comb {\NN} k$. Assume further that 
the indicator functions of $\X_{i_1, \ldots, i_k}$ belong to a compact subset $\mathcal K$ of $L^1(\Omega, \mu)$.  
Then, for any $\eps>0$ there exists an infinite set $J\subset \NN$ such that 
\[
\mu\left( \bigcap_{(i_1,\ldots,i_k)\in J^{[k]}} X_{i_1 \ldots i_k}\right)\ge \lambda -\eps.
\]
\et

\bp 
Consider first the case $k=1$. 
By compactness of $\mathcal K$, for all $\eps>0$ there exist an increasing sequence $\{i_n\}$ and a set $X_\infty\subset X$, with $\mu(X_\infty)\ge \lambda$, such that 
\[
\mu\left(X_\infty\Delta X_{i_n}\right) \le \frac{\eps}{2^n} \qquad \forall n\in\NN.
\]
As a consequence, letting $J:=\{i_n:\, n\in\NN\}$ we have 
\[
\mu\left(\bigcap_{n\in \NN}X_{i_n}\right)\ge
\mu\left(X_\infty\cap \bigcap_{n\in \NN}X_{i_n}\right)\ge 
\mu\left( X_\infty\right) - \sum_{n\in \NN}\mu\left(X_\infty\Delta X_{i_n}\right) \ge \lambda - \eps.
\]
For $k>1$, we apply Lemma \ref{lempie} with 
\begin{eqnarray*}
M &=& \mathcal K\subset L^1(\Omega,\mu)
\\
f(i_1, \ldots ,i_k) &=& \chi^{}_{X_{i_1 \ldots i_k}}\in L^1(\Omega,\mu). 
\end{eqnarray*}
In particular, recalling Remark \ref{remtopo}, for all $\eps>0$ 
there exist $J=\sigma(\NN)$, $X_\infty\subset \Omega$, and $X_{i_1 \ldots i_m}\subset X$, 
for all $(i_1, \ldots ,i_m)\in J^{[m]}$ with $1\le m<k$, such that  
$\mu(X_\infty)\ge\lambda$ and for all $(i_1, \ldots ,i_k)\in J^{[k]}$ it holds
\begin{eqnarray*}
\mu\left( X_\infty \Delta X_{i_1}\right) &\le& \frac{\eps}{2^{\sigma^{-1}(i_1)}}
\\
\mu\left( X_{i_1 \ldots i_m}\Delta X_{i_1 \ldots i_{m+1}}\right) &\le& \frac{\eps}{2^{\sigma^{-1}(i_{m+1})}}\,.
\end{eqnarray*}
Reasoning as above, it then follows
\begin{eqnarray*}
&&\mu \left(X_\infty\Delta \bigcap_{(i_1, \ldots ,i_k)\in J^{[k]}} X_{i_1 \ldots i_k}\right) \le 
\\
&&\sum_{i_1\in\NN} \mu\left(X_\infty\Delta X_{i_1}\right) 
+ \sum_{i_1<i_2} \mu\left( X_{i_1}\Delta X_{i_1 i_2}\right) +  
\\
&&\cdots  +  \sum_{i_1<\cdots<i_{k}} \mu\left( X_{i_1 \ldots i_{k-1}}\Delta X_{i_1 \ldots i_k}\right)
\le C(k)\eps\,,
\end{eqnarray*}
where $C(k)>0$ is a constant depending only on $k$. 
Therefore
\begin{eqnarray*}
\mu\left( \bigcap_{(i_1,\ldots,i_k)\in J^{[k]}} X_{i_1 \ldots i_k}\right) &\ge&
\mu\left( X_\infty\cap \bigcap_{(i_1,\ldots,i_k)\in J^{[k]}} X_{i_1 \ldots i_k}\right)
\\
&\ge& \mu\left( X_\infty\right) - \mu\left( X_\infty\Delta \bigcap_{(i_1,\ldots,i_k)\in J^{[k]}} X_{i_1 \ldots i_k}\right)
\\
&\ge&
\lambda - C(k)\eps.
\end{eqnarray*}
\ep

Notice that from Theorem \ref{taiut} it follows that 
Problem \ref{problemstrong} has a positive answer if there exist an infinite $J\subseteq\NN$
and sets $\widetilde \X_{i_1, \ldots, i_k}\subseteq X_{i_1 \ldots i_k}$
with $(i_1,\ldots,i_k)\in J^{[k]}$, such that 
$\mu\left(\widetilde \X_{i_1, \ldots, i_k}\right)\ge \lambda$ for some $\lambda>0$,
and the indicator functions of $\widetilde \X_{i_1, \ldots, i_k}$ belong to a compact subset of $L^1(\Omega,\mu)$.

\br\label{remper}
We recall that, when $\Omega$ is a compact subset of $\R^n$ 
and the perimeters of the sets $\X_{i_1, \ldots, i_k}$ are uniformly bounded,
then the family $\chi^{}_{\X_{i_1, \ldots, i_k}}$ has compact closure in $L^1(\Omega,\mu)$ 
(see for instance \cite[Thm. 3.23]{AFP:00}). In particular, 
if the sets $\X_{i_1, \ldots, i_k}$ have equibounded Cheeger constant, 
i.e. if there exists $C>0$ such that 
\[
\min_{E\subset \X_{i_1, \ldots, i_k}}\frac{{\rm Per}(E)}{|E|}\le C
\qquad \forall (i_1,\ldots ,i_k)\in\comb {\NN} k,
\]
then Problem \ref{problemstrong} has a positive answer.
\er


\appendix

\section{A topological Ramsey theorem} 

The following metric version of Ramsey theorem reduces to the classical Ramsey theorem when $M$ is finite. 

\bl \label{lempie} Let $M$ be a compact metric space, let $k\in \NN$, 
and let $f: \comb {\NN} k\to M$.  
Then there exists an infinite set $J\subset \NN$ such that the limit
\[
\lim_{\underset{(i_1, \ldots, i_k) \in \comb J k } {(i_1, \ldots, i_k) \to +\infty} } f(i_1, \ldots , i_k) 
\]
exists. 
\el

\bp Notice first that the thesis is trivial for $k=1$, since the space $M$ is compact. 
Assuming that the thesis holds for some $k\in \NN$, we want to prove it for $k+1$. So let $f\colon \comb {\NN} {k+1} \to M$. By inductive assumption, for all $j\in\NN$ there exist a infinite set $J_j\subset\NN$ and a point $x_j\in M$ such that $x_j=\lim_{i_1, \ldots, i_k \to\infty}f(j,i_1, \ldots , i_k)$, with $(i_1, \ldots , i_k)\in [J_j]^k$. 
Possibly extracting further subsequences we can also assume that 
\begin{equation}\label{eqJ}
d(x_j, f(j,i_1, \ldots , i_k)) \leq 1/2^j
\end{equation}
for all $(i_1, \ldots , i_k)\in \comb {J_j} k$. 
Moreover, by a recursive construction, we can 
assume that $J_{j+1}\subseteq J_j$. Now define $\tau \in \Ic$ by choosing $\tau(0) \in \NN$ and inductively $\tau(n+1)\in J_{\tau(n)}$. Since $J_{j+1}\subset J_j$ for all $j$, this implies $\tau(m)\in J_{\tau(n)}$ for all $m>n$. 
By compactness of $M$, there exists $\lambda \in \Ic$ and a point $x\in M$ such that $x_{\tau(\lambda(n))} \to x$ for $n\to \infty$. Take $J = \im (\tau \circ \lambda)$. The thesis
follows the triangle inequality $d(x, f(j,i_1, \ldots, i_k)) \leq d(x,x_j) + d(x_j, f(j,i_1,\ldots, i_k))$, noting that if $j<i_1<\ldots <i_k$ are in $J$, then $i_1, \ldots, i_k\in J_j$ (so Equation \ref{eqJ} applies).
\ep

Note that in Lemma \ref{lempie}, the condition $(i_1, \ldots, i_k) \to +\infty$ is equivalent to $i_1\to \infty$ (since $i_1<i_2<\ldots<i_k$). We would like to strengthen Lemma \ref{lempie} by requiring the existence of all the partial limits 
\[
x = \lim_{i_{j(1)} \to \infty} \lim_{i_{j(2)} \to \infty} \cdots \lim_{i_{j(r)} \to \infty} x_{i_1\ldots i_k} 
\]
where $1 \leq r\leq k$ and $(i_{j(1)}, \ldots, i_{j(r)})\in \comb J r$ is a subsequence of $(i_1, \ldots, i_k)\in \comb J k$. Note that the existence of all these $2^{k-1}$ partial limits does not follow from Lemma \ref{lempie}. For instance $\lim_{(i,j) \to \infty} \frac {(-1)^j}{i+1} = 0$ but $\lim_{i \to \infty} \lim_{j\to \infty} \frac {(-1)^j}{i+1}$ 
does not exist. 

To prove the desired strengthening it is convenient to introduce some terminology. 
Let $\overline\NN = \NN\cup \{\infty\}$ be the one-point compactification of $\NN$. 
Given a distance $\delta$ on $\NN$, we consider on $\comb {\NN} k$ the induced metric 
\[
\delta_k((n_1, \ldots, n_k), (m_1, \dots, m_k)) := \max_i \delta (n_i,m_i)\,.
\] 
Given $\sigma\in \Ic$, let $\sigma_* \colon \comb {\NN} k \to \comb {\NN} k$ be the induced map defined by $\sigma_* (n_1, \dots, n_k) := (\sigma(n_1), \ldots, \sigma (n_k))$. Given $f\colon \comb {\NN} k$, by the following theorem there is an infinite $J\subset \NN$ such that all the partial limits of $f \rest {\comb J k}$ exist. Moreover the arbitrarity of $\delta$ shows that we can impose an arbitrary modulus of convergence on all the partial limits of $f\circ \sigma_*$, where $\sigma \in \Ic$ is an increasing enumeration of $J$.

\bt \label{lempie2} Let $M$ be a compact metric space, let $k\in \NN$, 
and let $f: \comb {\NN} k\to M$. Then, for any distance $\delta$ on $\overline\NN$
there exists $\sigma\in\Ic$ such that 
$f\circ\sigma_*: \comb {\NN} k\to M$ is $1$-Lipschitz.
As a consequence, it can be extended to a $1$-Lipschitz function on the closure of $\comb {\NN} k$ in $\overline{\NN}^k$. 
\et

\bl \label{metric} Let $\delta$ be a metric on $\overline\NN$. Then there is another metric $\delta^*$ 
on $\overline\NN$ such that 
\begin{enumerate}
\item $\delta^*(x,y) \leq \delta(x,y)$ for all $x,y$. 
\item $\delta^*$ is \emph{monotone} in the following sense:  
$\delta^*(x',y') \leq \delta^*(x,y)$ for all $x,x, y, y'$, provided $x <\min(y,x',y')$.
\item $\eps^*(x)\ge\eps^*(y)$ for all $x\le y$, where 
\begin{equation}\label{reps}
\eps^*(x):=\min_{y\ge x+1}\delta^*(x,y).
\end{equation}
\end{enumerate}
\el

\bp 
We shall define a distance of the form $\delta^*(x,y) = \delta(\psi(x),\psi(y))$ for a suitable strictly increasing function $$\psi \colon \overline\NN\to \overline\NN\ .$$
To this aim, let us consider, for any $x\in \overline\NN$, the diameter of the interval $[x,\infty]\cap\overline\NN$ 
\begin{equation}\label{eps} 
\eta(x) := \max_{x\leq y\leq z} \delta (y,z), 
\end{equation}
and the point-set distance from $x$ to the interval $[x+1,\infty]\cap\overline\NN$
\begin{equation} \label{eta} 
\eps(x):=\min_{y\ge x+1} \delta(x,y).
\end{equation}
Since $\eps(x)>0$ for all $x<\infty$ and  $\eta(x)=o(1)$ as $x\to\infty$, there exists a recursively defined, strictly increasing function $\psi \colon \overline\NN\to \overline\NN$ such that for any $x\in\NN$

\begin{eqnarray}
\eta(\psi(x)) &\le &   \eps(x)    \\
\nonumber
 \eta(\psi(x+1)) &\le& \eps(\psi(x))      \,.  
\end{eqnarray}
As a consequence, the distance
\[
\delta^*(x,y) := \delta(\psi(x),\psi(y))
\]
verifies, for all $x<y\le\infty $
\[
\delta^*(x,y)=\delta(\psi(x),\psi(y)) \le \eta(\psi(x))\le \eps(x)\le \delta(x,y), 
\]
and, assuming also $x<x'\le\infty$ and $x<y'\le\infty$,
\begin{eqnarray*}
\delta^*(x',y') &= &\delta(\psi(x'),\psi(y')) \le \eta(\psi(x')) \le  \eta(\psi(x+1)) 
\\
\nonumber & \leq &\eps(\psi(x))\le \delta(\psi(x),\psi(y))=\delta^*(x,y)\,.
\end{eqnarray*}

To prove the last statement we observe that 
\[
\eps^*(x)\ge \eps(\psi(x))\ge \eta(\psi(x+1))\ge \eps^*(x+1).
\]
\ep 

\begin{proof}[Proof of Theorem \ref{lempie2}] 
By Lemma \ref{metric} we can assume that $\delta$ is monotone
in the sense of Lemma \ref{metric} (2).

We proceed by induction on $k$. 
When $k=1$, consider the function $\eps(n):=\min_{m\ge n+1}\delta(n,m)$ 
as in \eqref{reps}. By compactness of $M$ there exist $x\in M$ and a subsequence $f \circ \sigma$ 
of $f$ converging to $x$  with the property
\begin{equation}
d_M\left( f\left(\sigma n \right), x \right) \le \frac{\eps(n)}{2} \,.
\end{equation}
Recalling Lemma \ref{metric} (3), for $n \neq m$ we have 
\begin{equation}
d_M \left( f \left(\sigma n\right), f \left( \sigma m \right) \right) 
\leq  \frac{\eps(n)+\eps(m)}{2}\leq \delta (n,m)\,. 
\end{equation}
So $f\circ \sigma$ is 1-Lipschitz. 

Now assume inductively that the thesis holds for some $k\in \NN$, and let us prove it for $k+1$. So let $f\colon \comb {\NN} {k+1} \to M$. We need to prove the existence of $\sigma\in \Ic$ such that
\begin{equation}\label{lip} 
d_M\left( f\left(\sigma_*(n,\boldsymbol m)\right), 
f\left(\sigma_*(n',\boldsymbol m')\right)\right) \le \delta_{k+1}((n,\boldsymbol m), (n',\boldsymbol m'))
\end{equation}
for all $(n,\boldsymbol m) \in \comb {\NN} {k+1}$ and $(n',\boldsymbol m') \in \comb {\NN} {k+1}$, where $\boldsymbol m = (m_1, \ldots, m_k)$ and $\boldsymbol m' = (m'_1, \ldots, m'_k)$. 

Given $n\in \NN$ define $f_n\colon \comb {\NN} k\to M$ by
\begin{equation}
f_n(\boldsymbol m ) := \begin{cases} f(n,\boldsymbol m) & \text{if $n<m_1$,} 
\\
\perp & \text{if $n\geq m_1$}
\end{cases} 
\end{equation}
where $\perp$ is an arbitrary element of $M$. Note that the condition $n<m_1$ is equivalent to $(n,\boldsymbol m) \in \comb {\NN} {k+1}$. 

By inductive assumption, for all $n\in\NN$ there exists $\theta_n \in\Ic$ such that 
$f_n \circ {\theta_n}_* \colon \comb {\NN} k \to M$ is $1$-Lipschitz. By a recursive construction, we can also assume that $\theta_{n+1}$ is a subsequence of $\theta_n$, namely $\theta_{n+1} = \theta_n \circ \gamma_n$ for some $\gamma_n \in \Ic$. Indeed to obtain $\theta_{n+1}$ as desired it suffices to apply the induction hypothesis to $f_{n+1} \circ {\theta_n}_* \colon \comb {\NN} k\to M$ rather than directly to $f_{n+1}$. 

Since $f_n \circ {\theta_n}_*$ is 1-Lipschitz, there exist the limit 
\[
g(n) := \lim_{\min (\boldsymbol m)\to \infty} f(n, {\theta_n}_*(\boldsymbol m))
\]
Passing to a subsequence we can further assume that all the values of $f_n \circ \theta_n$ are 
within distance $\eps (n)/4$ from its limit, namely:
\begin{equation}\label{g} 
d_M \left( g \left(n \right),  f\left(n, \theta_n (\left(\boldsymbol m\right) \right) \right) < \frac{\eps (n)}{4}\,.
\end{equation} 
Let $J_n := \theta_n(\NN) \subset \NN$ and let $\tau\in \Ic$ be such that:
\begin{equation}\label{tau} 
\tau(n+1) \in J_{\tau(n)}
\end{equation}
It then follows that 
\begin{equation} \label{tau2} 
\forall n,m \in \tau(\NN)  \quad m>n \Longrightarrow m\in J_n \,. 
\end{equation} 
For later purposes we need to define $\tau(n+1)$ as an element of $J_{\tau(n)}$ bigger than its $n+1$-th element, namely $\tau(n+1) > \theta_{\tau(n)}(n+1)$. So, for the sake of concreteness, we define inductively $\tau(0):=0$ and $\tau(n+1) := \theta_{\tau(n)}(n+2)$. It then follows that: 
\begin{equation}\label{tau3} \forall i,j\in \tau(\NN) \; \forall k \in \NN \quad j>i, j\geq k \Longrightarrow
\tau(j)> \theta_{\tau(i)} (k) \,.
\end{equation} 
Reasoning as in the case $k=1$, there is $\lambda \in \Ic$ and $x_\infty\in M$ such that 
\begin{equation}\label{lambda}
d_M\left( g\left(\tau \left( \lambda \left( n \right) \right) \right), x_\infty \right) < \frac{\eps (n)}{4}
\end{equation} 

Now define $\sigma:= \tau \circ \lambda \in \Ic$. Note that $\sigma (\NN) \subset \tau (\NN)$ so (\ref{tau2}) and (\ref{tau3}) continue to hold with $\sigma$ instead of $\tau$. 
We claim that $f \circ \sigma_* \colon \comb {\NN} {k+1}\to M$ is 1-Lipschitz. 

As a first step we show that 
\begin{equation} \label{sigma} 
\exists \boldsymbol k> \boldsymbol m \st
(f \circ \sigma_*) (n,\boldsymbol m) = (f_{\sigma(n)} \circ \theta_{\sigma(n)})(n, \boldsymbol k)
\end{equation} 
where $\boldsymbol k> \boldsymbol m$ means that $k_i>m_i$ for all respective components. To prove (\ref{sigma}) recall that $(f\circ \sigma_*)(n, \boldsymbol m) = f(\sigma (n), \sigma(m_1), \ldots, \sigma (m_k))$. 
Since $n<\min (\boldsymbol m)$, by (\ref{tau2}) the elements $\sigma(m_1), \ldots, \sigma(m_k)$ are in the image of $\theta_{\sigma(n)}$, namely for each $i$ we have $\sigma(m_i) = \theta_{\sigma(n)} (k_i)$ for some $k_i \in \NN$. Moreover applying (\ref{tau3}) we must have $k_i > m_i$. The proof of (\ref{sigma}) is thus complete. 

It follows from (\ref{sigma}) and (\ref{g}) that $(f \circ \sigma_*) (n,\boldsymbol m)$ is within distance 
$\eps(\sigma(n))/4$ from its limit $g(\sigma(n))$, which in turn is within distance $\eps(n)/4$ 
from its limit $x_\infty$ by (\ref{lambda}). We thus proved: 
\begin{equation}
d_M\left(f \left( \sigma_* \left(n, \boldsymbol m \right) \right) , x_\infty \right) < 
\frac 1 4 \eps(\sigma(n)) + \frac 1 4 \eps(n) \,.
\end{equation}
Recalling that for $x\neq y$ we have $\eps(x)+\eps(y) \leq 2\delta(x,y)$, it follows that for $n\neq n'$ the left-hand side of (\ref{lip}) is bounded by 
$[\delta (\sigma(n), \sigma(n')) + \delta (n,n')]/2$, which in turn is $\leq \delta(n,n')$ by monotonicity of $\delta$. 

If remains to prove (\ref{lip}) in the case $n=n'$. Given $\boldsymbol m, \boldsymbol m'$ as in (\ref{lip}), we apply (\ref{sigma}) to get 
$\boldsymbol k > \boldsymbol m, \boldsymbol k' > \boldsymbol m'$ with 
$(f \circ \sigma_*) (n,\boldsymbol m) = (f_{\sigma(n)} \circ \theta_{\sigma(n)})(n, \boldsymbol k)$ and $(f \circ \sigma_*) (n,\boldsymbol m') = (f_{\sigma(n)} \circ \theta_{\sigma(n)})(n, \boldsymbol k')$.  

Using the monotonicity of $\delta$ and the fact that $f_{\sigma(n)} \circ \theta_{\sigma(n)}$ is 1-Lipschitz, it follows that:
\begin{equation}
d_M\left( f\left(\sigma_*(n,\boldsymbol m)\right), 
f\left(\sigma_*(n,\boldsymbol m')\right)\right) 
\le \delta_{k}(\boldsymbol k, \boldsymbol k') \le \delta_k (\boldsymbol m, \boldsymbol m')\,.
\end{equation}
\end{proof}

\br \label{remtopo}
 Theorem \ref{lempie2} implies that there exists an infinite set $J = \sigma(\NN)\subset\NN$ such that,
for all $0\le m< k$ and $(i_1,\ldots ,i_m)\in J^{[m]}$, there are limit points 
$x_{i_1\ldots i_m}\in M$ with the property
\[
x_{i_1\ldots i_m} = \lim_{\underset{(i_1\ldots i_k)\in J^{[k]}}
{(i_{m+1}, \ldots, i_k) \to\infty}}x_{i_1\ldots i_k},
\]
where we set $x_{i_1\ldots i_k} := f\left(i_1,\ldots ,i_k\right)$.
Moreover, by choosing the distance $\delta(n,m)=\eps |2^{-n}-2^{-m}|$,
we may also require  
\begin{equation*}
d_M\left(x_{i_1\ldots i_m},x_{i_1\ldots i_k}\right) \le \frac{\eps}{2^{\sigma^{-1}(i_{m+1})}}
\qquad \forall(i_1,\ldots,i_k)\in J^{[k]}.
\end{equation*}
\er


\section{Exchangeable measures} \label{app-exch}

Let $\Lambda$ be a compact metric space. 
We recall a classical notion of \emph{exchangeable measure} due to De Finetti \cite{Defin}, showing some equivalent conditions.

\bprop\label{proequiv}
Given $m\in \mathcal M^1(\Lambda^{\NN})$, the following conditions are equivalent:
\begin{enumerate}
\item[a)] $m$ is $\Sc$-invariant; 
\item[b)] $m$ is $\Ij$-invariant; 
\item[c)] $m$ is $\Ic$-invariant.
\end{enumerate}
\eprop

\bd
If $m$ satisfies one of these equivalent conditions we say that $m$ is \emph{exchangeable}.  
\ed 

Notice that an exchangeable measure is always shift-invariant, while there are  
shift-invariant measures which are not exchangeable.
To prove Proposition \ref{proequiv} we need some preliminary results concerning measures satisfying condition (c).  

\bd
Given $m\in \mathcal M(\Lambda^{\NN})$ and $f\in L^p(\Lambda^{\NN})$, with $p\in [1,+\infty]$, 
we let 
\[
\tilde f = E\left(f\vert \A_s\right)\in L^p(\Lambda^{\NN})
\]
be the conditional probability of $f$ with respect to the $\sigma$-algebra $\A_s$ 
of the shift-invariant Borel subsets of $\Lambda^{\NN}$.
In particular, $\tilde f$ is shift-invariant, and by Birkhoff's theorem (see for instance \cite{P:82}) we have 
\[
\tilde{f}= \lim_{n\to\infty} \frac{1}{n}\sum_{k=0}^{n-1}f\circ {{\sf S}^* }^k\,,
\]
where the limit holds almost everywhere and in the strong topology of $L^1(\Lambda^{\NN})$. 
\ed 

\bl \label{weak} 
Assume that $m\in {\ca M}^1(\Lambda^{\NN})$ is $\Ic$-invariant. Then  
for all $f\in L^\infty(\Lambda^{\NN}, m)$ we have
\begin{equation}\label{flim}
\tilde{f}= \lim_{n\to\infty} f\circ {{\sf S}^*}^n \,,
\end{equation}
where the limit is taken in the weak$^*$ topology of $L^\infty(\Lambda^{\NN})$, namely 
for every $g \in L^1(\Lambda^{\NN},m)$ we have 
\begin{equation}\label{eqlam}
\lim_{n\to \infty} \int_{\Lambda^{\NN}} g \, \big (f \circ {{\sf S}^*}^n\big) \,dm = \int_{\Lambda^{\NN}} g \tilde f  \, dm 
\end{equation}
\el 

\bp 
It suffices to prove that $\lim_{n\to\infty} f\circ {{\sf S}^*}^n$ exists, since in that case it is necessarily equal to the (weak$^*$) limit of the arithmetic means $\frac{1}{n}\sum_{k=0}^{n-1}f\circ {{\sf S}^* }^k$, and therefore to $\tilde f$ (since $\tilde f = \lim_{n\to \infty} \frac{1}{n}\sum_{k=0}^{n-1}f\circ {{\sf S}^* }^k$ in an even stronger topology). 
Since the sequence $f\circ {{\sf S}^*}^n$ is equibounded in $L^\infty(\Lambda^{\NN},m)$, 
it is enough to prove \eqref{eqlam}
for all $g$ in a dense subset $D$ of $L^1(\Lambda^{\NN})$. 
We can take $D$ to be the set of those functions $g \in L^1(\Lambda^{\NN},m)$ depending on finitely many coordinates (namely $g(x) = h(x_1,\ldots, x_r)$ for some $r\in \NN$ and some $h \in L^1(\Lambda^r, m)$). 
The convergence of \eqref{eqlam} for $g(x) = h(x_1,\ldots, x_r)$
follows at once from the fact that $\sigma \cdot m = m$ for all $\sigma \in \Ic$, 
which implies that the quantity in \eqref{eqlam} is constant for all $n>r$.
Indeed to prove that $\int_{\Lambda^{\NN}} g \, \big(f \circ {{\sf S}^*}^n\big) \,dm = \int_{\Lambda^{\NN}} g\, \big(f \circ {{\sf S}^*}^{n+l} \big) \,dm$ it suffices to consider the function $\sigma \in \Ic$ which fixes $0,\ldots, r-1$ and sends $i$ to $i+l$ for $i\geq r$.   
\ep

We are now ready to prove the equivalence of the conditions in the definition of exchangeable measure. 

\begin{proof}[Proof of Proposition \ref{proequiv}]
Since $\Sc\subset \Ij$ and $\Ic\subset \Ij$,
the implications ${\rm b)}\Rightarrow {\rm a)}$ and ${\rm b)}\Rightarrow {\rm c)}$ are obvious. The implication ${\rm a)}\Rightarrow {\rm b)}$ is also obvious since it is true on the Borel subsets of $\Lambda^{\NN}$ of the form $\{x \in \Lambda^{\NN} \st x_{i_1} \in A_1, \ldots, x_{i_r} \in A_r\}$, which generate the whole Borel $\sigma$-algebra of $\Lambda^{\NN}$.

Let $m\in {\ca M}^1(\Lambda^{\NN})$ be $\Ic$-invariant, and let us prove that $m$ is $\Ij$-invariant. So let $\sigma \in \Ij$. We must show that 
\begin{equation}\label{eqlim}
\int_{\Lambda^{\NN}} g \,dm = 
\int_{\Lambda^{\NN}} g \circ \sigma^* \,dm\,,
\end{equation}
for all $g\in C(\Lambda^{\NN})$. 
It suffices to prove \eqref{eqlim} for $g$ in a dense subset $D$ of $C(\Lambda^{\NN})$. So we can assume that $g(x)$ has the form $g_0(x_0)\cdot \ldots \cdot g_r(x_r)$ for some $r\in \NN$ and $g_1, \ldots, g_r \in C(\Lambda)$. Note that $g_i(x_i) = (g_i \circ P_i) (x)$ where $P_i \colon \Lambda^{\NN} \to \Lambda$ is the projection on the $i$-th coordinate. Since $P_i = P_0 \circ {\sf S}^* $ where ${\sf S}^* $ is the shift, we can apply Lemma \ref{weak} to obtain
\begin{equation*}
\int_{\Lambda^{\NN}} g \,dm = 
\int_{\Lambda^{\NN}} \widetilde{g_1\circ P_{1}}\cdots \widetilde{g_r\circ P_{1}}\,dm\,.
\end{equation*}
Reasoning in the same way for the function $g\circ \sigma^*$, we finally get 
\[
\int_{\Lambda^{\NN}} g\circ \sigma^* \,dm = 
\int_{\Lambda^{\NN}} \widetilde{g_1\circ P_{1}}\cdots \widetilde{g_r\circ P_{1}}\,dm = 
\int_{\Lambda^{\NN}} g \,dm\,.
\]
\end{proof}

\bd We say that $m \in {\ca M}^1(\Lambda^{\NN})$ is {\em asymptotically exchangeable} if 
the limit $m' = \lim_{\underset{\theta \in \Ic} {\min \theta \to \infty}} \theta \cdot m$ exists in ${\ca M}^1(\Lambda^{\NN})$ and is an exchangeable measure.  \ed

\br
Note that if $m$ is asymptotically exchangeable, then:
\begin{eqnarray}
m' & := & \lim_{\underset{\theta \in \Ic} {\min \theta \to \infty}} \theta \cdot m \\
 & = & \lim_{k\to \infty} {\sf S}^k \cdot m \,.
\end{eqnarray}
However it is possible that $\lim_{k\to \infty} {\sf S}^k \cdot m$ exists and is exchangeable, and yet $m$ is not asymptotically exchangeable. As an example one may start with the Bernoulli probability measure $\mu$ on $2^{\NN}$ with $\mu(\{x_i = 0\}) = 1/2$ and then consider the conditional probability $m(\cdot) = \mu(\cdot | A)$ where $A\subset 2^{\NN}$ is the set of those sequences $x\in 2^{\NN}$ satisfying $x_{(n+1)^2} = 1 - x_{n^2}$ for all $n$. 
\er

\br \label{j}
If $m$ is asymptotically exchangeable and $m'= \lim_{k\to \infty} {\sf S}^k \cdot m$, then
for all $r\in\NN$ and $g_1,\ldots,g_r\in C(\Lambda)$
we have
\begin{equation}\label{eqrico}
\lim_{\underset{(i_1, \ldots, i_r) \in \comb {\NN} r} {i_1\to +\infty} }
\int_{\Lambda^{\NN}} g_1(x_{i_1})\,\cdots\,g_r(x_{i_r}) \,dm
=\int_{\Lambda^{\NN}} g_1(x_{1})\,\cdots\,g_r(x_{r}) \,dm^\prime .
\end{equation}
\er

\bt \label{lemexch}
Given $m\in {\ca M}^1(\Lambda^{\NN})$ there is $\sigma \in \comb \omega \omega$ such that $\sigma \cdot m$ is asymptotically exchangeable. \et

\bp Fix $m\in {\ca M}^1(\Lambda^{\NN})$. Given $r \in \omega$ consider the function $f\colon \comb {\omega} r\to {\ca M}^1(\Lambda^r)$ sending $\iota$ to $\iota \cdot m\in {\ca M}^1(\Lambda^r)$. By Lemma \ref{lempie} there is an infinite set $J_r \subset \omega$ such that
\begin{equation} \label{limit}
\lim_{\underset{\iota \in \comb {J_r} r}{\min(\iota) \to \infty}} \iota \cdot m
\end{equation}
exists in ${\ca M}^1(\Lambda^r)$. By a diagonal argument we choose the same set $J = J_r$ for all $r$. Let $\sigma \in \Ic$ be such that $\sigma (\NN) = J$. We claim that $\sigma \cdot m$ is asymptotically exchangeable. To this aim consider $m_k := {\sf S}^k \cdot \sigma \cdot m \in {\ca M}^1(\Lambda^{\NN})$. By compactness there is an accumulation point $m'\in {\ca M}^1(\Lambda^{\NN})$ of $\{m_k\}_{k\in \NN}$. We claim that
\begin{equation} \label{mprime}
\lim_{\underset{\theta \in \comb {J} \omega}{\min(\theta) \to \infty}} \theta \cdot \sigma \cdot m = m'\, ,
\end{equation}
hence in particular $m_k \to m'$ (taking $\theta = {\sf S}^k$). Note that the claim also implies that $m'$ is exchangeable. Indeed, given an increasing function $\gamma \colon \NN \to \NN$, to show $\gamma \cdot m' = m'$ it suffices to replace $\theta$ with $\theta \circ \gamma$ in equation \eqref{mprime}. Since the subset of $C(\Lambda^{\NN})$ consising of the functions depending on finitely many coordinates is dense, it suffices to prove that for all $r\in \NN$ and $\iota \in \comb {\NN} r$ the limit
\begin{equation} \label{proj}
\lim_{\underset{\theta \in \comb {J} \omega}{\min(\theta) \to \infty}} \iota \cdot \theta \cdot \sigma \cdot m 
\end{equation}
exists in ${\ca M}^1(\Lambda^r)$ (the limit being necessarily $\iota \cdot m'$). This is however just a special case of equation \ref{limit}. \ep

 We give below some representation results for exchangeable measures. 
First note that if $\Lambda$ is countable, a measure $m \in {\ca M}^1(\Lambda^{\NN})$ is determined by the values it takes on the sets of the form $\{ x \st x_{i_1} = a_1, \ldots, x_{i_r} = a_r \}$.  

\bl \label{countable}
If $\Lambda$ is countable, a measure $m \in {\ca M}(\Lambda^{\NN})$ is
exchangeable if and only if it admits a representation of the following form. There is a probability space $(\Omega, \mu)$ (which in fact can be taken to be $(\Lambda^{\NN}, m)$) and a family $\{\psi_a\}_{a\in \Lambda}$ in $L^\infty(\Omega, \mu)$ such that 
for all $i_1< \ldots < i_r$ in $\NN$ we have 
\begin{equation}\label{repr} 
m\left( \left\{ x \st x_{i_1} = a_1, \ldots, x_{i_r} = a_r \right\} \right) =  \int_\Omega \psi_{a_1}\cdot \ldots \cdot \psi_{a_n} d\mu \, .
\end{equation}
\el

\bp 
Since the right-hand side of the equation does not depend on $i_1, \ldots, i_r$ a measure $m \in {\ca M}^1(\Lambda^{\NN})$ admitting the above representation is clearly exchangeable. Conversely if $m$ is exchangeable it suffices to take $\psi_a = \widetilde {\chi_a}$ where $\chi_a$ is the characteristic function of the set $\{x \st x_0 = a\}$. We can in fact obtain the desired result by a repeated application of Equation \eqref{eqlam} after observing that the characteristic function $\chi_{\{ x \st x_{i_1} = a_1, \ldots, x_{i_r} = a_r \}}$ is the product $\chi_{\{ x \st x_{i_1} = a_1\}}\cdot \ldots \cdot \chi_{\{x_{i_r} = a_r \}}$ and 
$\chi_{\{x \st x_i = a\}} =  \chi_{a} \circ ({\sf S}^* )^i$. 
\ep 

\bc \label{zero}
If $\Lambda$ is countable and $m \in {\ca M}^1(\Lambda^{\NN})$ is exchangeable, then $m(\{x \in \Lambda^{\NN} \st x_0 = x_1\}) \neq 0$. 
\ec 
\bp By \eqref{repr} $m (\{x \in \Lambda^{\NN} \st x_0 = x_1\}) = \sum_{a\in \Lambda} \int {\psi_a}^2 d\mu \neq 0$.  
\ep 

\bc \label{1p} 
If $p\in \NN$ and $m \in {\ca M}^1(p^{\NN})$ is exchangeable, then $m(\{x \in \Lambda^{\NN} \st x_0 = x_1\}) \geq \frac 1 p$. 
\ec 
\bp Write $m(\{x \in \Lambda^{\NN} \st x_0 = x_1\}) = \sum_{a\in \Lambda} \int_\Omega {\psi_a}^2$ and apply the Cauchy-Schwarz inequality to the linear operator $\sum \int$ on $p\times \Omega$ to obtain 
\begin{equation} 
\left(\sum_{a<p} \int_\Omega {\psi_a}^2 d\mu \right) \cdot \left(\sum_{a<p} \int_\Omega 1 d\mu \right)  \geq \left(\sum_{a<p} \int_\Omega \psi_a d\mu \right)^2
\end{equation}
which gives the desired result. 
\ep

Thanks to a theorem of De Finetti, suitably extended in \cite{HS:55}
there is an integral representation \emph{\`a la Choquet} for the exchangeable measures on $\Lambda^{\NN}$, where $\Lambda$ is a compact metric space. 
More precisely, in \cite{HS:55} it is shown that the extremal points of the (compact) convex set of all exchangeable measures 
are given by the product measures $\sigma^{\NN}$, with $\sigma \in \mathcal M^1(\Lambda)$. 
As a consequence, Choquet theorem \cite{Choquet} provides an integral representation 
for any exchangeable measure $m$ on $\Lambda^{\NN}$,  
i.e. there is a probability measure $\mu\in\mathcal M^1(\Lambda)$ such that
\begin{equation}\label{choquetgen}
m = \int_{\mathcal M^1(\Lambda)} 
\sigma^{\NN}\, d\mu(\sigma) \,.
\end{equation}
When $\Lambda$ is finite, i.e. $\Lambda= p= \{0,\ldots,p-1\}$
for some $p\in \NN$, we can identify ${\mathcal M}^1(\Lambda)$ with the 
symplex $\Sigma_p$ of all $\lambda\in [0,1]^p$ such that $\sum_{i=0}^{p-1}\lambda_i=1$. 
Given $\lambda\in \Sigma_p$, we denote by $B_\lambda$ the product  
measure on $p^{\NN}$, namely the unique measure making all the events $\{x\st x_i = a\}$  independent with measure $B_\lambda(\{x \st x_i=a\})=\lambda_a$. 
In this case, \eqref{choquetgen} becomes 
\begin{equation}\label{choquetbern}
m = \int_{\Sigma_p} B_{\lambda}\,d\mu(\lambda)\,,
\end{equation}
where $\mu$ is a probability measure on $\Sigma_p$.

We finish this excursus on exchangeable measures with the following result: 

\bprop \label{lemshift}
Let $m\in \mathcal M^1(\Lambda^{\NN})$ be exchangeable, then 
for all $f\in L^1(\Lambda^{\NN})$ the following conditions are equivalent:
\begin{enumerate}
\item[a)] $f$ is $\Sc$-invariant; 
\item[b)] $f$ is $\Ij$-invariant; 
\item[c)] $f$ is shift-invariant.
\end{enumerate}
\eprop

\bp
Since $\Sc\subset \Ij$ and $s\in \Ij$,
the implications ${\rm b)}\Rightarrow {\rm a)}$ 
and ${\rm b)}\Rightarrow {\rm c)}$ are obvious.

In order to prove that ${\rm a)}\Rightarrow {\rm b)}$, we let $\mathcal F = \{ \sigma\in\Ij:\,
f=f\circ \sigma^*\}$, which is a closed subset of $\Ij$ containing $\Sc$.
Then, it is enough to observe that 
$\Sc$ is a dense subset of $\Ij\subset \NN^{\NN}$, 
with respect to the product topology of $\NN^{\NN}$, 
so that $\mathcal F = \overline{\Sc}=\Ij$.

Let us prove that ${\rm c)}\Rightarrow {\rm a)}$. Let $\sigma\in \Sc$ and let $n$ be such that 
$\sigma(i)=i$ for all $i\ge n$. It follows that ${{\sf S}^*}^k\circ \sigma^* = {\sf S}^k$, for all $k\ge n$. 
As a consequence, for $m$-almost every $x\in \Lambda^{\NN}$ it holds
\[
f \circ \sigma^* (x) = f\circ {{\sf S}^*}^n \circ \sigma^* (x) = f \circ {{\sf S}^*}^n (x) = f(x),
\]
where the first equality holds since the measure $m$ is $\Sc$-invariant.
\ep

Notice that from Proposition \ref{lemshift} it follows that $\tilde f$ is $\Ij$-invariant 
for all $f\in L^1(\Lambda^{\NN})$.
In particular, for an exchangeable measure, the $\sigma$-algebra of the shift-invariant sets
coincides with the (a priori smaller) $\sigma$-algebra of the $\Ij$-invariant sets.

\end{document}